\documentclass[final]{autart}
\usepackage{natbib}

 \usepackage{parskip} 
 
 \usepackage{mathtools}
\usepackage{thmtools}
\usepackage{amsmath}
\usepackage{amssymb}
\usepackage[amsthm]{ntheorem}
\usepackage{enumerate}
\usepackage{graphicx}
 
\usepackage{verbatim}
\usepackage{url}
\usepackage{amsfonts} 
\usepackage{graphicx}
\usepackage{color}
 
\usepackage{epsfig}

 \newcommand{\Int}{\operatorname{int}}

 \newcommand{\diag}{\operatorname{diag}}
 
 \newcommand{\tr}{\operatorname{tr}}
 
\newcommand{\st}{\, | \,}
  \DeclareMathOperator{\vol}{vol}
 \DeclareMathOperator{\vspan}{span}

\newcommand\mycom[2]{\genfrac{}{}{0pt}{}{#1}{#2}}
\newcommand{\trace}{\operatorname{trace}}

\newcommand{\spanop}{\operatorname{span}}

\newcommand*\diff{\mathop{}\!\mathrm{d}}

 \newtheorem{Theorem}{Theorem}
 \newtheorem{Lemma}{Lemma}

 \newtheorem{Corollary}{Corollary}
  \newtheorem{Definition}{Definition}
 
 \newtheorem{Proposition}[Theorem]{Proposition}
 \newtheorem{Remark} {Remark}
 \newtheorem{Example} {Example}

\newcommand {\R}{\mathbb R}

\newcommand {\C}{\mathbb C}

\newcommand{\be}{\begin{equation}}
\newcommand{\ee}{\end{equation}}

\newcommand{\updt}[1]{{\color{black}#1}}

\begin{document}

\begin{frontmatter}
\title{$k$-Contraction:  Theory and Applications}
\thanks[footnoteinfo]{
This research  was partially supported by  a research grant  from  the Israel Science Foundation.}
 \author[First]{Chengshuai Wu},
 \author[First]{Ilya Kanevskiy},
 \author[Third]{Michael Margaliot}
\address[First] {C. Wu and I. Kanevskiy are with the School of Elec. Eng.,
		Tel Aviv University, Tel-Aviv 69978, Israel.}
\address[Third]{M. Margaliot (Corresponding Author) is with the Department of  Elec. Eng.-Systems and the Sagol School of Neuroscience,  Tel-Aviv
 University, Tel-Aviv 69978, Israel. E-mail:  \texttt{michaelm@tauex.tau.ac.il}}

\maketitle

\begin{abstract}
A dynamical system is called contractive if any two solutions approach one another at an exponential rate. 
More precisely, the dynamics contracts lines at an exponential rate.
This property implies highly ordered  asymptotic behavior including entrainment to 
time-varying periodic vector fields and, in particular, global asymptotic stability for time-invariant vector fields. 
Contraction theory  
has found numerous applications in systems and control theory because there exist 
 easy to verify sufficient conditions, based on matrix measures, guaranteeing contraction.

Here, we provide a geometric generalization of contraction theory called~\updt{$k$-contraction}.
A dynamical system is called~\updt{$k$-contractive} if the dynamics 
contracts  $k$-parallelotopes at an exponential rate.
For~$k=1$ this reduces to standard contraction. 

We describe easy to verify sufficient conditions for~\updt{$k$-contraction} based
on a   matrix measure of   the~$k$th additive compound of the Jacobian of the vector field. 
We also describe applications of  the seminal  work of Muldowney and Li,
that  can be interpreted in the framework of~\updt{$2$-contraction}, to systems and control theory.
\end{abstract}

 \begin{keyword}
Asymptotic stability, contraction analysis, matrix measures, variational equation, 
 entrainment, 
 compound matrices. 
  \end{keyword}

\end{frontmatter}

\section{Introduction}
  
    Contraction theory provides \updt{powerful}
    tools for analyzing the asymptotic behavior
		of time-varying nonlinear dynamical systems~\citep{LOHMILLER1998683, sontag_cotraction_tutorial, forni2014}. Unlike  Lyapunov methods, it studies the difference
		between any pair of solutions rather than   convergence to a specific   solution. 
		If the difference converges to zero then this implies highly ordered behavior. For example, if the state-space includes an equilibrium~$e$ then any solution is attracted to~$e$ implying global \updt{exponential} asymptotic stability. 
		More generally, if the vector field is time-varying and~$T$-periodic \updt{and the state-space is compact and convex} then there exists a unique~$T$-periodic solution~$\gamma$ and any solution  converges to~$\gamma$~\citep{entrain2011,LOHMILLER1998683}. 
		In other words, the system \emph{entrains} to the periodic excitation modeled by the time-varying  vector field.
		
		   Sufficient conditions for contraction can be derived   using
	a differential Lyapunov function~\citep{LOHMILLER1998683, forni2014}
	or by showing that some  matrix measure of the Jacobian of the vector field is uniformly negative~\citep{coppel1965stability,sontag_cotraction_tutorial}.  \updt{For a given matrix measure}, the latter condition is easy to check. 
		Contraction theory has found numerous applications in the field of 
		systems and control including:  control synthesis for regulation~\citep{pavlov2008incremental}
		and tracking~\citep{wu2019robust},  
		observer design \citep{doi:10.1002/aic.690460317,sanfelice2011convergence, aghannan2003intrinsic},   synchronization \citep{slotine2005study}, robotics \citep{manchester2018unifying}, multi-agents systems~\citep{russo2009solving}, and systems biology~\citep{RFM_entrain,entrain2011}.

	There is a large body of work on various
	generalizations of contraction theory. 
	Examples include contraction  with respect to~(w.r.t.) time- and space-dependent norms, that  
	are particularly relevant
    for systems whose trajectories evolve on manifolds~\citep{forni2014}.
	The recent paper by~\cite{jafarpour2020weak} considers contraction w.r.t. a 
	  seminorm. This is closely related to partial contraction 
		(or convergence to  an invariant 
		linear subspace)~\citep{slotine2005study}.  
    Another generalization is based on the fact that
     contraction guarantees strong asymptotic properties, like stability and entrainment, and thus it is often enough to consider systems that become contractive after some transient~\citep{cast_book,3gen_cont_automatica}.
    Another related line of work~\citep{forni_diff_diss} considers systems that are monotone w.r.t. ellipsoidal norms. These are not necessarily
		contractive systems, but the quadratic structure of the norm implies that they satisfy a form of partial
		contraction.

    In this paper, we present a geometric generalization   
		  called \emph{\updt{$k$-contraction}}. This is motivated by the seminal 
			work of~\cite{muldo1990} who used what we call here \updt{2-contraction} 
			to derive   generalizations of   results of Poincar\'e, Bendixson, and Dulac  on planar systems to
			higher-dimensional systems (see also~\cite{li1993bendixson, li1996geometric}). 
			The results of Muldowney and his colleagues proved very useful in analyzing 
			mathematical models for the spread of epidemics (see, e.g.,~\cite{SEIR_LI_MULD1995}). Indeed, these models typically include at least two equilibrium points corresponding to the disease-free and the endemic steady-states. Thus, they cannot be contractive.  However, they are sometimes \updt{$2$-contractive}  and this can be used to analyze their asymptotic behavior.

    To explain the notion of~\updt{$k$-contractive} systems in the simplest setting, consider a time-varying linear system. Fix~$k+1$ different initial conditions on the unit simplex, and an initial time~$t_0$.
		The corresponding solutions define at any time~$t\geq t_0$ a $k$-parallelotope.
		The system is called~\updt{$k$-contractive} if  the  volume 
		of this parallelotope decays to zero at an exponential rate. For~$k=1$ this reduces to standard contraction. For   nonlinear systems, \updt{$k$-contraction} is defined by considering a $k$-parallelotope on the \emph{tangent space}~\citep{carmo1992riemannian}.

	The tools needed to define and analyze~\updt{$k$-contraction} include the multiplicative and additive compound matrices~\citep{muldo1990}. The latter  also play  an important role in the theory of totally positive dynamical systems~(see the recent tutorial~\citep{margaliot2019revisiting} and also~\citep{Eyal_k_posi}). These notions  are not necessarily
	well-known in the systems and control community, and we try to provide  here  
	a self-contained exposition of  \updt{$k$-contraction} and its analysis using these tools. 
	
	To provide intuition, we begin with two simple linear examples. 
	The analysis of nonlinear systems is based on studying the associated variational equation which is
	a linear time-varying system. 
	\begin{Example}\label{exa:sipo}
	Consider the LTI system
	\be\label{eq:lti22}
	\dot x=Ax, \text{ with } A\in \mathbb{R}^{2\times 2}.
	\ee
	Let~$x(t,x_0)$ denote the solution of~\eqref{eq:lti22} at time~$t$ for
	the initial condition~$x(0)=x_0$. 
	Pick~$u,v \in\R^2$. \updt{Consider the 
	parallelogram~$  \{ r_1 x(t,u)+r_2 x(t,v) \st r_1,r_2 \in [0,1]\}$.}
The area of this parallelogram is $|s(t)|$, where
	\[
	s(t):=  \det(\begin{bmatrix} x(t,u) & x(t,v )\end{bmatrix} ) .
	\]
	This gives
	\begin{align*}
	s(t)&=  \det(\begin{bmatrix} \exp(At) u & \exp(At) v \end{bmatrix} )\\&=
	\det(\exp(At))    \det(\begin{bmatrix} u
	&  v \end{bmatrix} )  \\
	&= \det(\exp(At))  s(0).
	\end{align*}
	By the Abel-Jacobi-Liouville identity \citep{teschl2012ordinary}, 
$
\frac{d}{dt} \det(\exp(At))= \tr(A) \det( \exp(At) ),
$	 where~$\tr(A)$ denotes  the trace of~$A$, 
	so~$ \dot s(t)= \tr(A) s(t) $, and~$	s(t)
 = \exp(\tr(A) t )  s(0).
$
Summarizing, the area of the parallelogram spanned by the solutions~$x(t,u)$ and~$x(t,v)$ of the LTI 
system~\eqref{eq:lti22} decays to zero at an exponential rate if and only if~(iff) $\tr(A)<0$. We then say that the two-dimensional system~\eqref{eq:lti22} is~\updt{$2$-contractive}. The condition~$\tr(A)<0$ is weaker than that needed for standard contraction, namely, that~$A$ is Hurwitz~\citep{sontag_cotraction_tutorial}. On the other-hand,  standard contraction   implies that $\tr(A) <0$, i.e., \updt{$2$-contraction}.
\end{Example}

To generalize these notions to systems whose trajectories evolve on $\mathbb{R}^n$, with~$n >2$,   requires the use of \updt{multiplicative and additive  compound matrices}. The next example demonstrates this.
	
\begin{Example}
Consider the LTI  system
	\be\label{eq:lti23}
	\dot x=Ax, \text{ with } A \in  \R^{3\times 3}.
	\ee
	Pick~two arbitrary initial conditions $u,v  \in \R^3$. Let 
$
	M:=\begin{bmatrix}u & v\end{bmatrix}= \begin{bmatrix} u_1 & v_1\\ u_2 & v_2\\u_3 & v_3 \end{bmatrix}.
$
 Recall that the $2$nd multiplicative compound matrix of~$M$, denote~$M^{(2)}$, is the matrix of all~$2\times 2$ minors of~$M$ ordered lexicographically (as explained in Section \ref{subsec:compound} below).  Specifically, 
$	\updt{M^{(2)}}= \begin{bmatrix}
		u_1 v_2 - u_2 v_1 \\ u_1v_3-u_3v_1 \\ u_2 v_3-u_3v_2
	\end{bmatrix}$.
	The entries of this vector are (up to a minus sign) the entries of the 
	 cross product~$u\times v$. Thus,~$\left |\updt{M^{(2)}}  \right |_2=| u \times v|_2$, where~$|\cdot|_2$ is the~$L_2$ norm,
	and this implies that~$\left |\updt{M^{(2)}} \right |_2$ is
	  the area of the 
	parallelogram determined by~$u,v$. 
	
\updt{Thus, the area of the parallelogram generated by~$x(t, u)$ and~$x(t,v)$ is $|s(t)|_2$, where~$s(t) := \begin{bmatrix} x(t,u)& x(t,v) \end{bmatrix}^{(2)} $. Recall that the multiplicative compound  satisfies $(PQ)^{(2)}=P^{(2)}Q^{(2)}$ for any~$P \in \mathbb{R}^{n \times m},Q \in \mathbb{R}^{m \times k}$. Thus, we have} 
\begin{align}\label{eq:sder}
s(t)
	&= \begin{bmatrix} \exp(At) u& \exp(At ) v \end{bmatrix}^{(2)}   \nonumber\\
	&= (\exp(At )  M )^{(2)}   \nonumber\\
	&=(\exp(At))^{(2)}  M^{(2)} \nonumber \\
	&= (\exp(At ) ) ^{(2)}  s(0).
\end{align}		
It is useful to derive a differential equation for~$s(t)$. 
The $2$nd   \emph{additive compound} of a square matrix~$P$ is defined as
\begin{align*}
P^{[2]}&:= \frac{d}{d \varepsilon} (I+\varepsilon P)^{(2)} |_{\varepsilon=0}\\
&=  \frac{d}{d \varepsilon} (\exp(\varepsilon P))^{(2)} |_{\varepsilon=0}. 
\end{align*}
The term additive is due to the fact that~$(P+Q)^{[2]}=
P^{[2]}+ Q^{[2]}$ for any~$P,Q$ in $\mathbb{R}^{k \times k}$. 
Now~\eqref{eq:sder} gives~$
\dot s(t) = A^{[2]} s(t). 
$
 Thus, if~$\mu (\cdot)$ denotes the
 matrix measure induced by a  norm~$|\cdot |$, and~$\mu ( A^{[2]} )\leq-\eta<0$  then~$
|s(t)| \leq \exp(-\eta t) |s(0)|.
$
We then say that the
three-dimensional system~\eqref{eq:lti23} is~\updt{$2$-contractive}.
\end{Example}
Summarizing, \updt{$k$-contraction} is related to the contraction 
of the volume of $k$-parallelotopes under the dynamics. 

	The remainder of this paper is organized as follows.
	The next section reviews several notions and results that 
	are required to analyze $k$-contraction including
	compound matrices and \updt{the volume of parallelotopes}.
	Section~\ref{sec:koc} introduces the new notion of~\updt{$k$-contraction}.
	Section~\ref{sec:app} describes applications of \updt{$k$-contraction} to several  problems from systems and control theory. 
	
	\section{Preliminaries}
This section describes several notions that are used later on.
We begin by reviewing multiplicative and additive compound matrices,
and then the relation between \updt{the volume of $k$-parallelotopes}  and  the multiplicative compound.
We also review the spectral properties of  compound matrices, and describe  
the relation between the stability of the LTV system~$\dot x(t)=A(t)x(t)$
 and the stability of the associated $k$th compound system.

 For two integers~$i,j$, with~$ i\leq j$, let~$[i,j]:=\{i,i+1,\dots,j\}$. 
 Let~$Q_{k,n}$ denote the set 
of increasing sequences of~$k$ numbers from~$[1,n]$
ordered lexicographically.
For example,~$Q_{2,3} =
\{ (1,2), (1,3), (2,3) \}
$.
 \updt{The 
 lexicographic order~$\prec_{le}$ is defined as follows. If~$ a, b$ are
 two sequences in~$Q_{k,n}$, and $a_i$~$[b_i]$  is the $i$th element of $a$~$[b]$ then $a \prec_{le} b$ if $a_j < b_j$,  where $j = \min \{ i \in [1,k] \st  a_i \neq b_i\}$.}
With a slight abuse of notation we will sometimes treat   a     sequence in~$Q_{k,n}$ as a set.  

\subsection{Compound matrices}\label{subsec:compound}
Given a matrix~$A\in\R^{n\times m}$ and~$k\in[1,\min\{n,m\}]$, 
recall that a minor of order~$k$ of~$A$ is the determinant of some~$k \times k$
submatrix of~$A$.  
Consider the
$\binom{n}{k} \binom{m}{k}  $
 minors of  order~$k$ of~$A$. 
Each such minor is defined by a set of row indices~$\kappa_i \in Q_{k,n}$ and column indices~$\kappa_j\in 
Q_{k,m}$. This minor 
is denoted by~$A(\kappa_i|\kappa_j)$.
For example, for~$A=\begin{bmatrix} 4&5   \\ -1 &4   \\0&3 
\end{bmatrix}$,  we have
$
A(\{ 1,3\} |\{1,2\})=\det  \begin{bmatrix} 4&5\\0&3
\end{bmatrix}   =12.
$

The~$k$th \emph{multiplicative compound matrix} 
of~$A$, denoted~$A^{(k)}$, is the~$\binom{n}{k}\times  \binom{m}{k}$ matrix
that includes all  the minors of order~$k$ ordered lexicographically, \updt{that is, if~$\kappa_\ell$ is the~$\ell$th sequence in~$Q_{k,n}$ then the $ij$th entry of $A^{(k)}$ is $A(\kappa_i|\kappa_j)$.}
For example, for~$n = m =3$ and~$k=2$,  the matrix~$A^{(2)}$ is  
\[
\begin{bmatrix}
	                    A(\{1,2\}|\{1,2\}) & A(\{1,2\}|\{1,3\}) & A(\{1,2\}|\{2,3\})\\
						A(\{1,3\}|\{1,2\}) & A(\{1,3\}|\{1,3\}) & A(\{1,3\}|\{2,3\})\\
						A(\{2,3\}|\{1,2\}) & A(\{2,3\}|\{1,3\}) & A(\{2,3\}|\{2,3\}) 
\end{bmatrix}.
\] 
By definition,  $A^{(1)}=A$ 
and if~$A \in \mathbb{R}^{n \times n}$ then~$A^{(n)}=\det(A)$.
If~$D$ is an~$n\times n$ diagonal matrix, i.e.~$D=\diag(d_1,\dots,d_n)$ then 
\[
D^{(k)}=\diag(\prod_{i=1}^k d_i  , 
(\prod_{i=1}^{k-1} d_i )d_{k+1}, \dots, \prod_{i=n-k+1}^n d_i ).
\] 
 In particular, $I^{(k)}$ is the~$r\times r$ identity matrix, with~$r:=\binom{n}{k}$. 
 
The Cauchy-Binet formula  (see, e.g.,~\cite[Thm.~1.1.1]{total_book}) asserts   
that 
\be\label{eq:prodab}
(AB)^{(k)}=A^{(k)} B^{(k)} 
\ee
 for any $A \in \R^{n \times p}$, $B \in \R^{p \times m}$, and~$k \in [1,\min\{n,p,m\} ]$.
When~$n=p=m=k $ this becomes the familiar formula~$\det(AB)=\det(A)\det(B)$.
If~$A$ is square and non-singular then~\eqref{eq:prodab} 
  implies that~$ I^{(k)}=(AA^{-1})^{(k)}=A^{(k)} (A^{-1})^{(k)}$,
	so~$(A^{(k)})^{-1}=(A^{-1})^{(k)}   $. 

\updt{Note that   any entry 
in~$A^{(k)} $ is a polynomial in  the entries of~$A$. } The $k$th \emph{additive compound matrix} of~$A \in \R^{n \times n}$
is   defined  by~$
				A^{[k]}:= \frac{d}{d \varepsilon}  (I+\varepsilon A)^{(k)} |_{\varepsilon=0}.
$
This implies that~$A^{[1]}=A$, and that 
\be\label{eq:poyrt}
(I+\varepsilon A)^{(k)}= I+\varepsilon   A^{[k]} +o(\varepsilon )  ,
\ee 
that is,~$\varepsilon  A^{[k]}$ is the first-order term in the Taylor series of~$(I+\varepsilon A)^{(k)}$.
\begin{Example}\label{exa:sdig}
If~$D=\diag(d_1,\dots,d_n)$ then
\[
(I+\varepsilon D)^{(k)}=\diag( \prod_{i=1}^k (1+\varepsilon d_i) 
 ,\dots,\prod_{i=n-k+1}^n (1+\varepsilon d_i)      ),
\]
  so~\eqref{eq:poyrt} gives 
$
D^{[k]}=\diag( \sum_{i=1}^k d_i 
 ,\dots,  \sum_{i=n-k+1}^n d_i
 )     .
$
\end{Example}
\begin{Example}
Consider the case~$n=3$ and~$k=2$. Then
\begin{align*}
						&(I+\varepsilon A)^{(2)}  =\begin{bmatrix}
						1+\varepsilon a_{11} & \varepsilon a_{12} & \varepsilon a_{13}\\ 
						\varepsilon a_{21} &1+\varepsilon a_{22} & \varepsilon a_{23}\\
						\varepsilon a_{31}&\varepsilon a_{32} & 1+\varepsilon a_{33}
						\end{bmatrix}^{(2)}\\
						&=I +\varepsilon  \begin{bmatrix}
						 a_{11}+a_{22}  &  a_{23} & -  a_{13}\\
						  a_{32} &  a_{11}+a_{33}  &   a_{12}\\
						- a_{31}&  a_{21} &  a_{22}+a_{33} 
						\end{bmatrix}+o(\varepsilon),
\end{align*}
so 
\begin{align*} 
							A^{[2]}&=
						\frac{d}{d \varepsilon}  (I+\varepsilon A)^{(2)} |_{\varepsilon=0}\nonumber \\&=
						\begin{bmatrix}
						 a_{11}+a_{22} &   a_{23} & -  a_{13}\\
						  a_{32} &  a_{11}+a_{33} &  a_{12}\\
						-  a_{31}&  a_{21} &  a_{22}+a_{33}
						\end{bmatrix}.
\end{align*}
\end{Example}

It follows from~\eqref{eq:poyrt} and the properties of the multiplicative compound that
\begin{align*}
I+\varepsilon(A+B)^{[k]}+o(\varepsilon) &=(    I+\varepsilon (A+B)  )  ^{(k)}\\
&=(I+\varepsilon  A)^{(k)}(I+\varepsilon B)^{(k)}+o(\varepsilon)\\
&=(I+\varepsilon A ^{[k]}+o(\varepsilon) ) (I+\varepsilon B^{[k]}+o(\varepsilon) ) \\
&=I+\varepsilon (A ^{[k]}  + B ^{[k]})+ o(\varepsilon  ),
\end{align*}
so
\updt{taking~$\varepsilon\to 0$ and using the continuity w.r.t.~$\varepsilon$} yields
\[
(A+B)^{[k]}= A^{[k]}+B^{[k]},
\]
thus justifying the term additive compound.

The matrix~$A^{[k]}$ can be described
 explicitly in terms of the entries~$a_{ij}$ of~$A$. 
\begin{Lemma}[\cite{schwarz1970, fiedler_book}]\label{lem:poltr}  
Let~$A\in\R^{n\times n}$. Fix~$k\in[1,n]$. The entry of~$A^{[k]}$ corresponding to~$(\kappa_i|\kappa_j)=(i_1,\dots,i_k|j_1,\dots,j_k) $  is:
\begin{itemize}
\item $\sum_{\ell=1}^k a_{i_\ell i_\ell}$ if~$i_\ell=j_\ell$ for all~$\ell\in[1,k]$;
\item $(-1)^{\ell+m} a_{i_\ell j_m}$ 	 if all the indices in~$\kappa_i$ and~$\kappa_j$ coincide except for a single index~$i_\ell\not = j_m $; and
\item $0$, otherwise. 
\end{itemize}
\end{Lemma} 
 
The first case in Lemma~\ref{lem:poltr}
   corresponds to diagonal entries  of~$A^{[k]}$.  
All the other entries of~$A^{[k]}$ are
 either zero or     an entry of~$A$
multiplied  by either plus or minus one. 
Two special cases of Lemma~\ref{lem:poltr} are:
\[
A^{[1]}=A \text{ and } A^{[n]}=\tr(A).
\]

\begin{Example}
 Consider the case~$n=4$, i.e.,~$A=\{a_{ij}\}_{i,j=1}^4$. 
Then Lemma~\ref{lem:poltr} yields
\begin{align*}
A^{[2]}&=\left [ \begin{smallmatrix} 
																 a_{11}+a_{22}   & a_{23} & a_{24} & -a_{13} & -a_{14} & 0 \\
																a_{32}      &a_{11}+a_{33} & a_{34} & a_{12} & 0 & -a_{14} \\   
																a_{42}     &a_{43}   &a_{11}+a_{44} & 0 & a_{12} & a_{13} \\   
																-a_{31}      &a_{21} & 0 & a_{22}+a_{33} & a_{34} &  -a_{24} \\
																-a_{41}      &0 & a_{21} & a_{43} & a_{22}+a_{44} &  a_{23} \\ 
																0      &-a_{41} & a_{31} & -a_{42} & a_{32} &  a_{33}+a_{44}  
  \end{smallmatrix}\right] , 
	\end{align*}
and
\begin{align}\label{eq:a3fop}
A^{[3]}&=\left[ \begin{smallmatrix}     
																a_{11}+a_{22}+a_{33} & a_{34} & -a_{24} & a_{14} \\
																a_{43}& a_{11}+a_{22}+a_{44} & a_{23} & -a_{13}  \\
																-a_{42}& a_{32} & a_{11}+a_{33}+a_{44} & a_{12}  \\
                                a_{41}& -a_{31}& a_{21}& a_{22}+a_{33}+a_{44} 
 \end{smallmatrix} \right]. 
\end{align}
The entry in the first row and third column of~$A^{[3]}$
corresponds to~$(\kappa_i|\kappa_j)=( \{1,2,3\}| \{1,3,4\})$,
 and since~$\kappa_i$ and~$\kappa_j$ coincide except for the entry~${i_2}=2$   
and~${j_3}=4$, this entry is~$(-1)^{2+3} a_{i_2 j_3} = -a_{24}$.
It is useful to index compound matrices using~$\kappa_i,
\kappa_j$.
For example,  we write
$
A^{[3]}( \{1,2,3\}| \{1,3,4\})= -a_{24}.
$
\end{Example}

\updt{ We next review a   ``duality relation'' of the additive compound. 
This is a  formula relating the 
 matrices~$A^{[k]}$ and~$A^{[n-k]}$ that both 
have dimensions~$r\times r$, where~$r:=\binom{n}{k}$. 
Let~$U_r\in\R^{r\times r}$ be the matrix with entries
\be\label{eq:defuij}
u_{ij}=\begin{cases}
(-1)^{j+1} , &\text{if } i+j=r+1,\\
0,& \text{ otherwise}.
\end{cases}
\ee
For example,
$
U_3=\begin{bmatrix}
0 & 0 &1\\
0 &-1 & 0\\
1& 0&0
\end{bmatrix}.
$
Note that~$U_r^T=U_r^{-1}$. \\

\begin{Proposition}\citep{mol_appl_comp} \label{prop:dual_nnk}
Let~$A\in\R^{n\times n}$.
Fix~$k\in[1,n-1]$, and let~$r:=\binom{n}{k}$. Then
\be\label{eq:dual}
(A^{[k]})^T+U_r^T A^{[n-k]}   U_r=\tr(A) I_r .
\ee
\end{Proposition}
For example, suppose that~$n=3$,  and~$A=\diag(\lambda_1,\lambda_2,\lambda_3)$. Then for~$k=2$, \eqref{eq:dual} becomes
\begin{align*}
     \diag(\lambda_1+\lambda_2, &\lambda_1+\lambda_3,\lambda_2+\lambda_3) 
    +\diag( \lambda_3,\lambda_2,\lambda_1 ) \\&=(\lambda_1+\lambda_2+\lambda_3) I_3.
\end{align*}

}

\updt{Our next goal is to provide a clear geometric interpretation for   
the  multiplicative compound of a matrix.}

\subsection{\updt{The volume of  a parallelotope}}
\updt{Fix~$k\in [1,n]$. The parallelotope generated by the vectors~$x^1,\dots,x^k \in \R^n$ is the set
\[
P(x^1, \dots, x^k):=\{ \sum_{i=1}^k r_i x^i \st r_i \in [0,1] \}.
\]
Note that~$P$ always includes the origin, and that~$P$ is the image of the unit~$k$-cube
under the matrix
 \[
 X:=\begin{bmatrix}
 x^1&\dots& x^k
 \end{bmatrix}\in\R^{n\times k} .
 \]

The \emph{Gram matrix} (see e.g. \cite[p. 441]{matrx_ana}) associated with~$x^1,\dots,x^k$
is the $k\times k$ symmetric matrix:
\begin{align}\label{eq:defgram}
    G ( x^1,\dots,x^k): =X^T X .
\end{align}
Note that $G$ is positive semi-definite, and is positive definite iff the vectors~$x^1, \dots, x^k$ are  linearly independent. 
For example, for~$k=2$ we have
\[ 
    G ( x^1, x^2): =
    \begin{bmatrix}
    |x^1|_2^2 & (x^1) ^Tx^2   \\
     (x^2) ^T x^1 & |x^2|_2^2   
    \end{bmatrix}, 
\]
where~$|\cdot|_2$ denotes the $L_2$ norm, i.e., the \emph{Euclidean norm}.

The volume of~$P(x^1,\dots,x^k)$,
denoted~$\vol(P(x^1,\dots,x^k))$, is defined in  a recursive manner.
\begin{Definition}\label{def:vol_parall}
 For~$k=1$,~$\vol(P(x^1)):=|x^1|_2$. For any~$k>1$,
\be\label{eq:def_vol}
\vol(P(x^1,\dots,x^k)):=\vol(P(x^1_,\dots, x^{k-1}) ) h,
\ee
where $h  \geq 0$ is the Euclidean distance from~$x^k$ to the subspace $\vspan\{x^1, \dots, x^{k-1}\}$.
\end{Definition}

Note that~$P(x^1_,\dots, x^{k-1})$ can be viewed as the~$(k-1)$-dimensional  ``base'' of~$P(x^1,\dots,x^ k)$, and $h$ is the related ``height''. Hence, Definition~\ref{def:vol_parall} has a clear geometric interpretation, and it reduces to the standard notions
of length, area, and volme if $k= 1, 2, 3$, respectively. 

The next result provides a simple
algebraic expression for the volume in terms of the Gram matrix. 
\begin{Proposition}
 \cite[Chapter~IX]{Gantmacher_vol1}
\label{prop:vol_via_gram}
The volume of~$P(x^1,\dots,x^k)$ satisfies 
\be\label{eq:volpeqn}
\vol(P(x^1,\dots,x^k))=\sqrt{\det(G(x^1,\dots,x^k))}. 
\ee
\end{Proposition}
} 

\updt{Note that in the special case where~$k=n$, the matrix~$X$ is a square matrix, and~\eqref{eq:volpeqn}
gives
\begin{align*}
(\vol(P(x^1,\dots,x^k)) ) ^2&= 
    \det(G(x^1,\dots,x^n)) \\
    &=  \det      \left (   X^T X
    \right ) \\
    &=(\det(X))^2,
\end{align*}
i.e. the well-known formula
    \[
    \vol(P(x^1,\dots,x^n))=|\det(\begin{bmatrix}
    x^1&\dots & x^n
    \end{bmatrix})|.
    \]
 
 To relate the volume of~$P(x^1,\dots,x^k)$ to the multiplicative compound, note that combining~\eqref{eq:defgram} and the Cauchy-Binet formula yields
 \begin{align*}
     \det(G(x^1,\dots,x^k))&=
     \det(X^TX)\\
     &=(X^T X) ^{(k)}\\
     &=(X ^{(k)})^T X ^{(k)}.
 \end{align*}
Since~$X\in\R^{n\times k}$, 
 $X^{(k)}$ is an~$\binom{n}{k}$ column vector, so~$\det(G(x^1,\dots,x^k))= |X ^{(k)}|_2^2$.
 Combining this with Prop.~\ref{prop:vol_via_gram} yields
the elegant formula 
\be\label{eq:vol_via_compound}
\vol(P(x^1,\dots,x^k))=|X ^{(k)}|_2. 
 \ee
When the vectors~$x^i$ depend on time and we consider asymptotic properties (e.g. convergence to zero of the volume as time goes to infinity), it is possible to use any vector norm~$|X ^{(k)}|$, as all norms on~$\R^n$ are equivalent.  }

\begin{figure}[t]
 	\begin{center}
		\includegraphics[width=8cm]{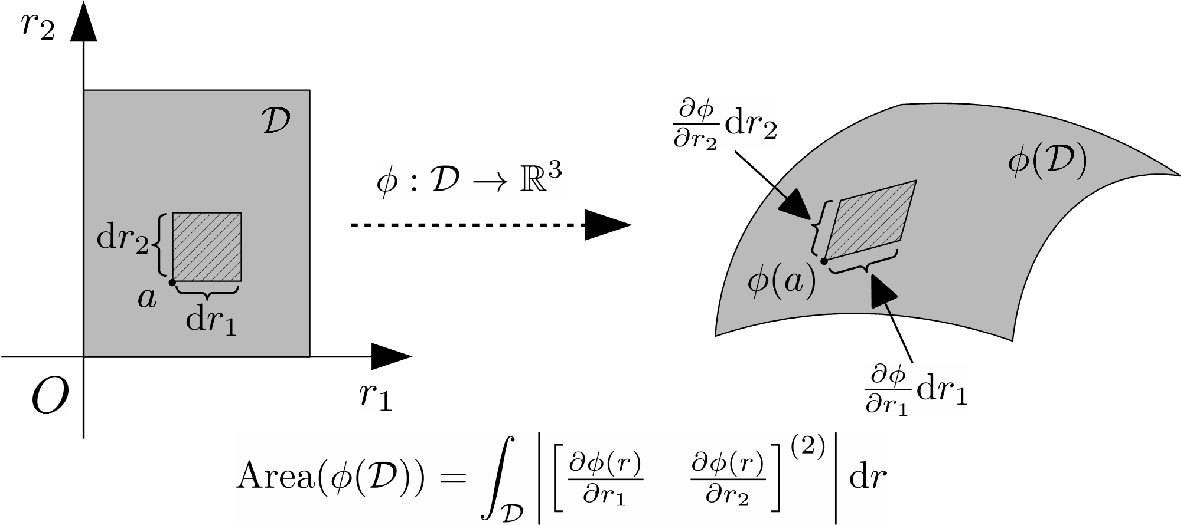} 
		 \caption{Using the multiplicative compound to compute the area of a parametrized body
		defined by~$\phi:\mathcal D\to\R^3$, with~$\mathcal D$ a rectangle in~$\R^2$.}
 		\label{wedge} 
 	\end{center}
 \end{figure}

\updt{Based on the above analysis}, the \updt{multiplicative compound} can be used to compute the volume of parameterized 
 bodies~\citep{muldo1990}. 
Consider a compact set~$\mathcal{D}  \subset \mathbb{R}^k$
 and a continuously  differentiable mapping~$\phi: \mathcal{D}  \rightarrow \mathbb{R}^n$, with $n\geq k$. This induces the set
\begin{equation}
\phi(\mathcal{D})  := \{    \phi (r) \st  r \in \mathcal{D} \}\subseteq \R^n.
\end{equation}
Since $\mathcal{D} $ is compact and $\phi (\cdot)$ is continuous, $\phi(\mathcal{D}) $ is closed. 
The volume of $\phi(\mathcal{D}) $ is  
\begin{equation}\label{eq:wjskcont}
\updt{\vol}( \phi(\mathcal{D})  ) = 
\int_{\mathcal{D} } \left \vert \updt{ \begin{bmatrix} \frac{\partial \phi (r)}{\partial r^1} & \cdots & \frac{\partial \phi (r)}{\partial r^k}   \end{bmatrix}^{(k)} } \right\vert \mathrm{d} r,
\end{equation}
(where we assume that the integral exists).
Fig.~\ref{wedge} illustrates this formula
for the case~$k = 2$ and~$n=3$.

\subsection{Spectral properties of compound matrices}

\updt{Recall that any~$A \in \R^{n \times n}$ admits a Jordan canonical form~\citep{achieser1992theory}: there exist~$T, J \in \C^{n \times n}$, with $T$ non-singular, such that~$A = T J T^{-1}$, where~$J$ has a Jordan normal form, and  in particular  is an upper-triangular matrix.
The diagonal entries of $J$, denoted $\lambda_i$, $i = 1, \dots, n$, are the eigenvalues of $A$.
Since~$A^{(k)} = T^{(k)} J^{(k)} (T^{(k)})^{-1}$   and~$J^{(k)}$ is upper triangular~\citep{muldo1990},
the diagonal entries of $J^{(k)}$ are   $\lambda_{i_1}\lambda_{i_2}\cdots \lambda_{i_k}$,
with~$1\leq i_1<\dots<i_k\leq n$,
and these are the eigenvalues of~$A^{(k)}$.

\begin{Example}
Suppose that~$Q\in\R^{n\times n}$ is symmetric and positive-definite.
Fix~$k \in[1,n]$. By definition,~$(Q^{(k)})^T =(Q^T)^{(k)}  $, so~$Q^{(k)}$ is symmetric. Every eigenvalue of~$Q$ is real and positive, and since every eigenvalue of~$Q^{(k)}$ is a product of~$k$ eigenvalues of~$Q$, every eigenvalue of~$Q^{(k)}$ is real and positive. We conclude that~$Q^{(k)}$ is positive-definite. 
\end{Example}

Let $u^i$ denote the eigenvector  of~$A$ corresponding to~$\lambda_i$, and let~$U: = \begin{bmatrix} u^{i_1} &  u^{i_2} & \cdots & u^{i_k} \end{bmatrix}$.
Then $A U = U \diag(\lambda_{i_1}, \lambda_{i_2}, \dots, \lambda_{i_k})$.   If $u^{i_1}, u^{i_2}, \dots, u^{i_k}$ are linearly  independent  then~$U^{(k)}$ is a nonzero vector and $A^{(k)} U^{(k)} = \prod_{j = 1}^k \lambda_{i_j} U^{(k)}$, i.e.  $U^{(k)}$ is an eigenvector of~$A^{(k)}$ corresponding to the eigenvalue~$\prod_{j = 1}^k \lambda_{i_j}$.

Similarly, $A = T J T^{-1}$ implies that $I + \varepsilon A = T (I + \varepsilon J ) T^{-1}$, and thus $(I + \varepsilon A)^{(k)} = T^{(k)} (I + \varepsilon J )^{(k)} (T^{(k)})^{(-1)}$. By \eqref{eq:poyrt},   $A^{[k]} = T^{(k)} J^{[k]} (T^{(k)})^{(-1)}$. 
Since $J $ is upper-triangular, so is~$J^{[k]}$, 
and the diagonal entries of $J^{[k]}$ are of the form~$\lambda_{i_1} + \lambda_{i_2}+ \cdots + \lambda_{i_k}$ according to Lemma~\ref{lem:poltr}.
Hence, every eigenvalue of~$A^{[k]}$ is the sum  of~$k$ eigenvalues of~$A$. }


 The standard tool for verifying   contraction and, as we will see below
also~\updt{$k$-contraction}, is matrix measures~\cite[Ch.~2]{coppel1965stability} 
(also called logarithmic norms~\citep{storm1975}). 

\subsection{Matrix measures}\label{subsec:matrixm}	 
Consider a vector norm $|\cdot|:\R^n\to\R_+$.
The \emph{induced matrix norm} $||\cdot||:\R^{n\times n}\to\R_+$ is 
$
||A|| := \max_{|x|=1} |Ax| $, and the
 induced  \emph{matrix measure} $\mu(\cdot):\R^{n\times n}\to\R $   is
\[
\mu(A) := \lim_{\varepsilon \downarrow 0} \frac{||I + \varepsilon A|| - 1}{\varepsilon} .
\]

Denote the~$L_1$,~$L_2$, and~$L_\infty$ vector norms by
$|x|_1:=\sum_{i=1}^n|x_i|$, 
$|x|_2:=\sqrt{\sum_{i=1}^n x_i^2}$, 
and~$|x|_\infty:=\max_{   i    } |x_i|$.  
The corresponding matrix measures are~\citep{vid}:
\begin{equation}\label{eq:matirxm} \begin{aligned} 
\mu_1(A) &= \max_{j}\left (a_{jj}+ \sum_{\substack{i=1 \\ i \neq j}}^n |a_{ij}| \right ) , \\
\mu_2(A) &= \lambda_1\left( \frac{A + A^T}{2} \right) , \\
\mu_{\infty}(A) &= \max_{i}\left (a_{ii}+ \sum_{\substack{j=1 \\ j \neq i}}^n |a_{ij}|\right ) , 
\end{aligned} \end{equation}
where~$\lambda_i (S)$ denotes the $i$-th largest eigenvalue of the symmetric matrix~$S$, that is,
\be\label{eq:prdeeifg}
\lambda_1(S) \geq \lambda_2(S) \geq \cdots \geq \lambda_n(S).
\ee

The 
 matrix measures for~$A^{[k]}$ are then~\citep{muldo1990}:
\begin{align}\label{eq:matirxm_k}
\mu_1(A^{[k]}) &= 
 \max_{(i)}\left ( \sum_{p=1}^k a_{i_p,i_p}   + \sum_{\substack{j \notin (i) }}(|a_{j,i_1}| + \cdots + |a_{j,i_k}|)\right ) ,\nonumber \\
\mu_2(A^{[k]}) &= \sum_{i=1}^k \lambda_i\left( \frac{A + A^T}{2} \right) , \\
\mu_{\infty}(A^{[k]}) &=
\max_{(i)}\left(  \sum_{p=1}^k a_{i_p,i_p} + \sum_{\substack{j \notin (i) }}(|a_{i_1,j}| + \cdots + |a_{i_k,j}|)\right) ,   \nonumber
\end{align}
where the maximum is taken over all $k$-tuples $(i) := \{i_1, \cdots, i_k \}\in Q_{k,n}$. Note that for~$k=1$, \eqref{eq:matirxm_k} reduces to~\eqref{eq:matirxm}.

The next subsection reviews several
 applications of compound matrices in dynamical systems described by~ODEs.
\subsection{Compound matrices and ODEs}
Consider the LTV system:
\be\label{eq:ltv}
\dot x(t)=A(t)x(t),
\ee
where~$A(t)$ is a continuous matrix function.
Then~$
x(t)=\Phi(t,t_0)x(t_0),
$
where~$\Phi$ is the transition 
 matrix corresponding to~\eqref{eq:ltv},  
satisfying 
\begin{align}\label{eq:ptrma}
						\frac{d}{dt} \Phi(t,t_0) =A(t) \Phi(t,t_0),\quad \Phi(t_0,t_0)=I. 
\end{align}
For the sake of simplicity, we always
assume    that the initial time is~$t_0=0$ and write~$\Phi(t)$ for~$\Phi(t,0)$. 

It is useful  to know how~$\Phi^{(k)}(t):=(\Phi (t))^{(k)}$
evolves in time.  Note that~$\Phi^{(k)}:\R_+\to\R^{r\times r}$, with~$r:=\binom{n}{k}$.
For example, 
\cite{schwarz1970} considered the following
 question: what  conditions  on~$A(t)$ guarantee that
 \emph{every} minor of~$\Phi(t)$ will be positive for all~$t>0$? 
In other words, $\Phi(t)$ is a \emph{totally positive matrix}~\citep{total_book} for all~$t>0$. 
When this holds~\eqref{eq:ltv} is called a \emph{totally positive differential system}~(TPDS). 
 Of course, the positivity of every minor of~$\Phi(t)$  is equivalent to the positivity of every entry in
 each of the matrices~$\Phi^{(1)}(t), \Phi^{(2)}(t),\dots,\Phi^{(n)}(t)$.

The additive compound arises naturally 
 when studying the dynamics of the multiplicative compound~$\Phi^{(k)}(t)$. Indeed,
for any~$\varepsilon>0$, 
\begin{align*}
      \Phi^{(k)} (t+\varepsilon)&=  (\Phi(t)+\varepsilon A(t)\Phi(t)   ) ^{(k)}+o(\varepsilon)\\
			&=(I+\varepsilon A(t)   ) ^{(k)} \Phi  ^{(k)}(t)+o(\varepsilon).
\end{align*}
Combining this with~\eqref{eq:poyrt} and the fact that~$\Phi(0)=I$ gives
\be\label{eq:povt}
\frac{d}{dt} \Phi^{(k)}(t)=A^{[k]}(t) \Phi^{(k)}(t), \quad  \Phi^{(k)}( 0)=I_r, 
\ee
where~$ A^{[k]}(t):=(A(t)) ^{[k]}$. 
In other words,   all the minors of order~$k$ of~$\Phi(t)$, stacked in the matrix~$\Phi^{(k)}(t)$, also follow
 a linear dynamics  
  with the matrix~$A^{[k]}(t) $. 
	Eq.~\eqref{eq:povt} is sometimes called the \emph{$k$th compound equation}
	(see e.g.~\cite{LI2000295}).

For a constant matrix~$A$, $\Phi(t)=\exp(At)$, so~\eqref{eq:povt} gives
\be\label{eq:gkpo}
	\exp(A^{[k]}t) =(\exp (At))^{(k)}. 
\ee
For~$k=n$,  Lemma~\ref{lem:poltr}  shows that~$A^{[n]}= \tr( A)$, whereas  
$(\exp (At))^{(n)}$ is the matrix that contains all the~$n\times n$ minors of~$\exp (At)$, that is,~$\det( \exp (At))$. Thus,~\eqref{eq:gkpo} generalizes  the Abel-Jacobi-Liouville identity.

It is useful to know how the $k$th compound equation~\eqref{eq:povt} changes under a coordinate transformation of~\eqref{eq:ltv}. Let~$T\in\R^{n\times n}$  be non-singular. Then
\be\begin{aligned} \label{eq:cor}
								(TAT^{-1})^{[k]}& =\frac{d}{d\varepsilon} (I+\varepsilon TAT^{-1})^{(k)}|_{\varepsilon=0}\\
								&=\frac{d}{d\varepsilon} (T (I+ \varepsilon A )T^{-1} )^{(k)}|_{\varepsilon=0}\\
&=T^{(k)}  A^{[k]} (T^{-1})^{(k)}.
\end{aligned}\ee

In the context of systems and control theory, it is important to understand the 
connections between stability of an  LTV system and of its associated~$k$th
 compound equation.

\subsection{Stability of an LTV system and of its $k$th compound equation}
As shown in~\cite[Corollary 3.2]{muldo1990}, 
under a certain boundedness  assumption there is an interesting 
relation between 
the stability of~\eqref{eq:ltv} and the stability of~$\dot y(t)=A^{[k]}(t)y(t)$. 
We state this result in a slightly modified form.

\begin{Proposition} \citep{muldo1990} 
\label{prop:sub}  
	Suppose that the LTV  system~\eqref{eq:ltv} is uniformly stable. Fix~$k\in[1,n]$.  
	Then 
	the following two conditions are equivalent. 
	\begin{enumerate}[(a)]
	\item \label{cond:nk1} The LTV system~\eqref{eq:ltv}
	admits an~$(n-k+1)$-dimensional linear subspace 
	$\mathcal{X} \subseteq \R^n$ such that
	\begin{equation} \label{eq:a0x}
	\lim_{t \to \infty }x(t, x_0) = 0 \text{ for any } x_0 \in \mathcal{X} .
	\end{equation}
	\item \label{cond:sysmkop} Every solution of  
	\begin{equation} \label{eq:kinsys}
	\dot{y}(t) = A^{[k]}(t)y(t)  
	\end{equation}
	converges to the origin.
\end{enumerate}
\end{Proposition}

For the sake of  completeness, we include  the  proof  in the  Appendix. 

 Prop.~\ref{prop:sub}  
  implies in particular that if every solution of~\eqref{eq:kinsys}
converges to the origin then for any~$\ell >k$ every solution of
$
	\dot{y}(t) = A^{[\ell]}(t)y(t)  
$
also converges to the origin.

\begin{Example}
Consider the simplest case, namely, the LTI system~$\dot x=Ax$, with~$A$ diagonalizable, that is, 
there exists a nonsingular matrix~$T$ such that~$TAT^{-1}=\diag(\lambda_1,\dots,\lambda_n)$. 
Assume that the real part of every~$\lambda_i$ is not positive, so that all the solutions are bounded. If there exist~$n-k+1$ eigenvalues with a negative real part (and thus~$k-1$  eigenvalues with a zero real part)
then: (1)~the dynamics admits an~$(n-k+1)$-dimensional linear subspace~$\mathcal{X}$ such that~$\lim_{t\to\infty} x(t,a)=0$ for any~$a\in \mathcal{X}$; and (2)~the sum of any~$k$ eigenvalues of~$A$ has a negative real part, so~$\dot y(t)=A^{[k]}y(t)$ is asymptotically stable.   
\end{Example}

\begin{Example}
	Consider the LTV system~$\dot{x}(t) = A(t)x(t)$ with~$n=2$ and
$
	A(t) = \begin{bmatrix}
	 -1 & 0 \\
	-\cos( t) &0  
	\end{bmatrix}.
$
	For any~$a\in\R^2$ the solution of this system is~$x(t,a)=\Phi(t)a$, with
	\begin{align*}
										\Phi(t)= \begin{bmatrix} \exp(-t)&0 \\ 
										               (-1+\exp(-t)(\cos(t)-\sin(t))) /2 &1\end{bmatrix}.
	\end{align*}
	This implies that the system is uniformly stable and that
	\be\label{eq:retyp}
	\lim_{t\to\infty}x(t,a)=\begin{bmatrix} 0 \\ 
										                  a_2-(a_1/2) \end{bmatrix}  .
	\ee
	The system is not contractive  w.r.t. any norm, as   not all solutions converge
	to the equilibrium~$0$. However,~$ A^{[2]}(t) =\tr(A(t)) \equiv -1$ (implying as we will see below that the
	system is~\updt{$2$-contractive}). In particular, for~$k=2$ Condition~\eqref{cond:sysmkop} in Prop.~\ref{prop:sub}
	holds. 
	By~\eqref{eq:retyp}, Condition~\eqref{cond:nk1}  also holds for the~$1$-dimensional linear subspace~$\mathcal{X}:=\spanop(\begin{bmatrix}2&1\end{bmatrix}^T)$. 
\end{Example}
  

We are now ready to  introduce  the main notion studied in this paper. 
\section{\updt{$k$-contraction}}\label{sec:koc}

Consider the time-varying  nonlinear system:
\be\label{eq:nonlinsys}
\dot x(t) = f(t,x),
\ee
where $f :\R_+\times \R^n \to \R^n$. 
We assume throughout  that the solutions evolve on a closed and convex  
state-space~$\Omega\subseteq\R^n$, and that for any initial condition~$a\in\Omega$, a unique solution~$x(t,a)$ exists and satisfies~$x(t,a)\in\Omega$  for all~$t\geq 0$. 
We also assume that~$f$ is continuously differentiable w.r.t. its second variable, 
and let~$J(t,x):=\frac{\partial}{\partial x}f(t,x)$ denote the Jacobian of~$f(t,x)$.

Pick~$a,b\in\Omega$.
Let~$h:[0,1]\to\Omega $ be the line~$h(r):=r a +(1-r)b$.
 Note that the convexity of~$\Omega$ implies that~$h(r)\in\Omega$ for all~$r\in[0,1]$. 
 Let~$
w(t,r):=\frac{\partial }{\partial r }x(t,h(r)).
$
Note that~$w(0,r)=\frac{\partial }{\partial r }x(0,h(r))= a-b$.
Intuitively speaking,~$w(t,r)$ measures 
 how a small change in the initial condition 
along the line~$h(r)$ affects the solution of~\eqref{eq:nonlinsys} at time~$t$. 
Then
\begin{align}\label{eq:pofuf}
\dot w(t,r)& :=\frac{d}{dt} w(t,r) \nonumber\\
& = \frac{\partial }{\partial r } \frac{d}{dt}  x(t,h(r))\nonumber \\
           &=  \frac{\partial }{\partial r } f(t, x(t,h(r) )) \nonumber  \\
					 &=   J(t, x(t,h(r)) ) w(t,r).
\end{align}					
This LTV system is the \emph{variational equation} associated with~\eqref{eq:nonlinsys}
along~$x(t,h(r))$,
as it describes how the variation between two initial conditions evolves  with time. 

If there exists a matrix measure such that
\be\label{eq:mudp}
\mu(J(t,z))\leq-\eta \text{ for all } t\geq 0 \text { and all } z\in\Omega
\ee
then it is not difficult to show~\citep{entrain2011}
 using~\eqref{eq:pofuf} that					
\[
			|x(t,a)-x(t,b)|\leq \exp(-\eta t) |a-b| \text{ for all } t\geq 0.
\]
If~$\eta> 0$ then this implies  contraction.

Our goal is to generalize these ideas  in the case where~\eqref{eq:mudp} is replaced 
by the more general condition~$
\mu( (  J(t,z)  ) ^{[k]})\leq -\eta  
$
for some~$k\in[1,n]$.
 It turns out that this  condition 
has a clear geometrical interpretation.
 To explain this, we first consider an 
LTV system  and   then proceed to explain 
the implications for the nonlinear system~\eqref{eq:nonlinsys}. 

\subsection{Linear time-varying systems}
We begin by considering the LTV system:
\be\label{eq:ltvw}
					\dot w(t)=A(t) w(t).
\ee
For the sake of simplicity, assume throughout 
that~$ A(t )$ is continuous in~$t$, but the extension to the
case of measurable and locally essentially bounded matrix functions is straightforward. 
This case is relevant, for example, when the dynamics depends on a control input.

\begin{Definition}\label{def:conlin}
Pick~$k\in[1,n]$. 
We say that~\eqref{eq:ltvw} is \updt{\emph{$k$-contractive}} if there exist~$\eta > 0$ and a vector norm $|\cdot|$ such that for
any~$a^1, \dotsm, a^k \in \R^n$, \updt{ the mapping $W(\cdot):\R_+\to  \R^{n \times k}$ defined by~$
W(t) : = \begin{bmatrix} w(t, a^1) & \dots & w(t, a^k)
\end{bmatrix}
$
satisfies
\be\label{eq:contdef}
| W^{(k)}(t)|
 \leq \exp(- \eta t)
 | W^{(k)}(0) |, \text{ for all }t \geq 0.
\ee}
\end{Definition}

\updt{In other words, under the dynamics  the volume of  any  $k$-parallelotope    decays to zero at an exponential rate. }
\begin{Example}
Consider~\eqref{eq:ltvw} with~$n=2$ and the constant matrix~$A=\begin{bmatrix}3& 0 \\0 &-4 \end{bmatrix}  $. Pick~$p,q\in\R^2$.
Then
\updt{
\begin{align*}
\vert \begin{bmatrix} w(t,p) & w(t,q) \end{bmatrix}^{(2)} & \vert
=\left |\begin{bmatrix} \exp(3t) p_1  & \exp(3t) q_1 \\  
\exp(-4t) p_2  & \exp(-4t) q_2  \end{bmatrix}^{(2)} \right | \\
= & \left |  \det \left ( \begin{bmatrix} \exp(3t) p_1  &  \exp(3t) q_1  \\ 
	\exp(-4t) p_2& \exp(-4t) q_2  \end{bmatrix}  \right) \right |\\
= & \exp(-t)\left |\begin{bmatrix} p & q \end{bmatrix}^{(2)}\right |,
\end{align*}}
so the system is~\updt{$2$-contractive} with~$\eta=1$. 
More generally, Example~\ref{exa:sipo} shows that when~$n=2$ and~$A\in\R^{2\times 2}$  is  a constant matrix then~\eqref{eq:ltvw} is~\updt{$2$-contractive} iff~$\tr(A)<0$. 
\end{Example}

An important advantage of standard contraction is that it admits  an easy to check
 sufficient condition based on matrix measures. 
The next result provides an easy   to check sufficient condition for \updt{$k$-contraction} of~\eqref{eq:ltvw} in terms of~$A^{[k]}(t )$. 

\begin{Theorem}\label{prop:kcontract}
If there exist $\eta > 0$ and a vector norm $|\cdot|$, with induced matrix measure $\mu : \R^{n \times n} \to \R$, such that
\be\label{eq:contrcond}
\mu(A^{[k]}(t )) \leq -\eta \text{ for all }  t \geq 0  
\ee
  then~\eqref{eq:ltvw} is \updt{$k$-contractive}. 
\end{Theorem}

\begin{proof}
For $k=1$  the definition of~\updt{$k$-contraction}
 reduces to standard contraction, and 
 condition~\eqref{eq:contrcond} reduces to the standard matrix measure sufficient condition for contraction, as~$A^{[1]}=A$. 
Consider the case $k>1$.  Pick $a^1, \dots, a^k \in \mathbb{R}^n$. Then
$   \dot W(t)  = A(t) W(t)$.
Hence, 
$
     \frac{d}{dt} W^{(k)}(t) = A^{[k]}(t ) W^{(k)}(t).
$
Now, using standard results on contraction,~\eqref{eq:contrcond} implies that
\[
|W^{(k)}(t)| \leq \exp(- \eta t) |W^{(k)}(0)| \text{ for all } t\geq 0, 
\]
and this completes the proof.
\end{proof}

\begin{Remark}
Consider~$A\in\R^{n\times n}$. There exists some matrix measure~$\mu$ such that~$\mu(A)<0$ iff~$A$ is Hurwitz~\citep{sontag_cotraction_tutorial}. 
Hence, there exists some matrix measure~$\mu$ such that~$\mu(A^{[k]})<0$ iff~$A^{[k]}$ is Hurwitz, that is,
iff the sum of every~$k$ eigenvalues of~$A$ has a negative real part. 
\end{Remark}

The next simple example  describes an LTI system that is ``on the verge'' of being~$2$-contractive.
\begin{Example}
Consider $\dot w = A w$ with 
$
A = \begin{pmatrix} 0 & 1 \\ -1 & 0 \end{pmatrix}.
$
Since $\frac{d}{dt} (x_1^2(t)  + x_2^2(t) ) = 0$, the solution for any~$x(0)$ is a circle with radius   
$
r := \sqrt{x_1^2 (t)  + x_2^2(t) } \equiv
 \sqrt{x_1 ^2 (0) + x_2 ^2(0)} . 
$
Let $W(t) := \begin{bmatrix} w(t, a^1) & w(t, a^2) \end{bmatrix}$. Then, 
\begin{align*}
    \frac{d}{dt}W^{(2)}(t) &= A^{[2]} W^{(2)}(t)
    = tr(A)W^{(2)}(t) 
    = 0.
\end{align*}
This agrees with the fact that the  area  of
the parallelotope generated by~$x(t,a^1)$ and~$x(t,a^2)$
  remains constant under the flow.
\end{Example}

\updt{We considered condition~\eqref{eq:contrcond} in the context of the~$n\times k$ matrix~$W(t)$ and the vector~$W^{(k)}(t)$. Yet, it also has implications for the~$n\times n$ transition matrix of the LTV~\eqref{eq:ltvw}. 
Let~$|\cdot|$ denote the vector norm that induces the matrix measure~$\mu$ in~\eqref{eq:contrcond} , and let~$||\cdot||$ denote the matrix norm induced by~$|\cdot|$.
\begin{Proposition}\label{prop:ulbound}
Let~$\Phi(t)$ denote the transition matrix corresponding to~\eqref{eq:ltvw}, that is,
\[
\dot \Phi (t)=A(t)\Phi(t),\quad \Phi(0)=I. 
\]
Fix~$k\in[1,n]$.
Then
\begin{align*}
    \exp(-\int_0^t \mu(-A^{[k]}(s)) &\diff s)     \leq    ||\Phi^{(k)}(t) ||  \\&  \leq \exp(\int_0^t \mu(A^{[k]}(s)) \diff s)   
\end{align*}
for all~$t\geq0$. 
In particular, 
if~\eqref{eq:contrcond} holds then
\[
||\Phi^{(k)}(t)|| \leq \exp(-\eta t),\text{ for all } t\geq 0 .
\]
\end{Proposition}
}
\updt{The proof follows by combining 
the fact that~$\dot \Phi^{(k)}(t) = A^{[k]}(t) \Phi^{(k)}(t)$,~$\Phi^{(k)}(0)=I$,
with Coppel's inequality.
 
For the particular case of the~$L_2$ norm 
we have~$||A||_2=\sqrt{\lambda_{\max}(A^TA) }$
and~$\mu_2(A)=(1/2)\lambda_{\max}(A+A^T)$, so~Prop.~\ref{prop:ulbound}
gives
\begin{align*}
   \exp  & (\int_0^t  \sum_{i=n-k+1}^n \lambda_i(A(s)+A^T(s)     ) \diff s)\\
   &\leq     \lambda_{\max} ( (\Phi^{(k)}(t))^T \Phi^{(k)}(t)   )\\&  \leq \exp(\int_0^t  \sum_{i=1}^k \lambda_i(A(s)+A^T(s)     ) \diff s)  . 
\end{align*}
This inequality has
been used  by~\cite{smith_Hausdorff_dim}
to bound  the Hausdorff dimension of chaotic attractors of nonlinear dynamical systems.
}

\updt{
The next result provides an inequality
relating    matrix measures of~$A^{[k]}$
and~$A^{[n-k]}$. Recall that if~$|\cdot|:\R^n\to\R_+$ is a vector norm,~$\mu:\R^{n\times n}$ is the induced matrix measure,  and
$P\in\R^{n\times n}$ is non-singular then the vector norm defined by~$|x|_P:=|Px|$ induces the matrix measure~$\mu_P(A)=\mu(PAP^{-1}) $.

\begin{Proposition}\label{prop:dual_vec_meas}
Let~$A\in\R^{n\times n}$.
Fix~$k\in [1,n-1 ] $, and let~$r:=\binom{n}{k}$. 
Let~$U_r$ be  the matrix defined in~\eqref{eq:defuij}.
Then for any matrix measure~$\mu$, we have 
\be\label{eq:inemp}
\mu((A^{[k]})^T)+\mu_{U_r^{T}}( A^{[n-k]} )
\geq \trace(A)  .
\ee
%
\end{Proposition}

\begin{proof}
Applying~$\mu$ on both sides  of~\eqref{eq:dual} and using  the   subadditivity   of the matrix measure yields~\eqref{eq:inemp}.
\end{proof}
}

 The next result shows that the sufficient condition 
for contraction w.r.t. some~$L_p$ norm, with~$p\in\{1,2,\infty\}$,
induces a ``graded structure''.
\begin{Corollary}\label{thm:graded}
If there exist $\eta > 0$ and~$p\in\{1,2,\infty\}$ 
such that
\be\label{eq:gradd}
\mu_p(A^{[k]}(t ) )\leq -\eta  \text{ for all }  t \geq 0  
\ee
then~\eqref{eq:ltvw} is~\updt{$\ell$-contractive} w.r.t. the~$L_p$ norm
for any~$\ell\geq k$. 
\end{Corollary}

\begin{proof}
We will prove this for the case~$p=2$. The proof for the cases~$p=1$ and~$p=\infty$ is based on similar arguments. 
Fix~$t\geq 0$ and let~$A=A(t)$. For~$p=2$   condition~\eqref{eq:gradd} is~$\sum_{i=1}^k \lambda_i(S)\leq-\eta<0$, 
where~$S:=(A+A^T)/2$ and the eigenvalues are ordered as in~\eqref{eq:prdeeifg}. 
This implies that~$\lambda_k(S)<0$ and thus~$\lambda_j(S)<0$ for any~$j \geq k$. 
Hence, for any~$\ell > k$ we have
\begin{align*}
\mu_2(A^{[\ell]})&=\sum_{i=1}^\ell \lambda_i(S)\\
&= \sum_{i=1}^k \lambda_i(S) + \sum_{i=k+1}^\ell \lambda_i(S) \\
&< \mu_2(A^{[k]}) \\
&\leq-\eta ,
\end{align*}  
and Theorem~\ref{prop:kcontract} implies that the system is~\updt{$\ell$-contractive} w.r.t. the~$L_2$ norm. 
\end{proof}

Theorem~\ref{prop:kcontract} can be used to provide new sufficient conditions for \updt{$k$-contraction}.
 The next two  results demonstrate this.
\begin{Proposition}
Suppose that $D(t) $ is diagonal and that there exists $k \in [1,n]$ such that the sum of every~$k$ diagonal entries of~$D(t)$ is smaller or equal to~$-\eta<0$ for all~$t\geq 0$. Then~$\dot x = D(t)x$ is \updt{$k$-contractive} w.r.t. the $L_s$ norm for any~$s\in\{1,2,\infty\}$.
\end{Proposition}
The proof follows from the fact that
\[
D^{[k]} = \begin{bmatrix} 
d_{11}+\dotsm +d_{kk} & 0 & \dotsm & 0 \\
\vdots & \ddots & & \vdots\\
 0 & \dotsm & & d_{pp}+\dotsm+ d_{nn}
\end{bmatrix},
\]
where $p:=n-k+1$. Thus, $\mu_1(D^{[k]}(t))\leq -\eta$ for all~$t$,
 and since~$D^{[k]} $ is diagonal,
$\mu_1(D^{[k]} ) =\mu_2(D^{[k]} )=\mu_\infty(D^{[k]} )$.

We can also derive a simple sufficient condition for \updt{$(n-1)$-contraction} in an~$n$-dimensional system. This is based on the following fact.
For~$M\in\R^{n\times n}$, let~$\tilde M$ denote the matrix with entries
\[
\tilde m_{ij}:=(-1)^{i+j} m_{n+1-i,n+1-j},\quad i,j \in[1,n].
\]
 \cite{schwarz1970} proved that if~$A\in\R^{n\times n}$ then
\be\label{eq:snm1}
			A^{[n-1]} = \tilde{B},
\ee
where~$B:=\tr(A) I -A^T$. For example, for $A \in \mathbb{R}^{4 \times 4}$, we have
$
B = 
\left[ \begin{smallmatrix}     
s - a_{11} & -a_{21} &- a_{31} &- a_{41}\\
-a_{12} & s - a_{22} &- a_{32} &- a_{42} \\
-a_{13} & -a_{23} & s - a_{33} &- a_{43} \\
-a_{14} & -a_{24} & -a_{34} &   s - a_{44}  \\
\end{smallmatrix} \right] ,
$
where $s : = \sum_{i=1}^4 a_{ii}$, so
\[
A^{[3]}  = \tilde{B} = 
\left[ \begin{smallmatrix}     
a_{11}+a_{22}+a_{33} & a_{34} & -a_{24} & a_{14} \\
a_{43}& a_{11}+a_{22}+a_{44} & a_{23} & -a_{13}  \\
-a_{42}& a_{32} & a_{11}+a_{33}+a_{44} & a_{12}  \\
a_{41}& -a_{31}& a_{21}& a_{22}+a_{33}+a_{44} 
\end{smallmatrix} \right] 
\]
 (compare with~\eqref{eq:a3fop}). 

\begin{Proposition}
Suppose that $A:\R_+\to   \R^{n \times n}$
satisfies    
\be\label{eq:tpyu}
      \sum_{\mycom{i=1}{i\not =\ell }}^{n}  \left( | a_{i \ell}(t)|+    a_{ii} (t)\right) \leq -\eta <0 ,
\ee
for all~$\ell \in[1,n]$ and all~$t\geq 0$.
Then~\eqref{eq:ltv} is~\updt{$(n-1)$-contractive} w.r.t. the~$L_\infty$ norm. 
\end{Proposition}

To show this, note that~\eqref{eq:snm1} implies 
that the sum of the entries of every
 row of~$A^{[n-1]}(t)$, with off-diagonal terms taken with absolute value, 
 is the expression on the left-hand side of~\eqref{eq:tpyu} for some~$\ell$,
so~$\mu_\infty(A^{[n-1]}(t))\leq -\eta$ for all~$t\geq 0$.

 We now turn to consider \updt{$k$-contraction} in nonlinear dynamical  systems.  
 
  \begin{figure*}
	\begin{center}
		\centering
	\includegraphics[scale=0.36]{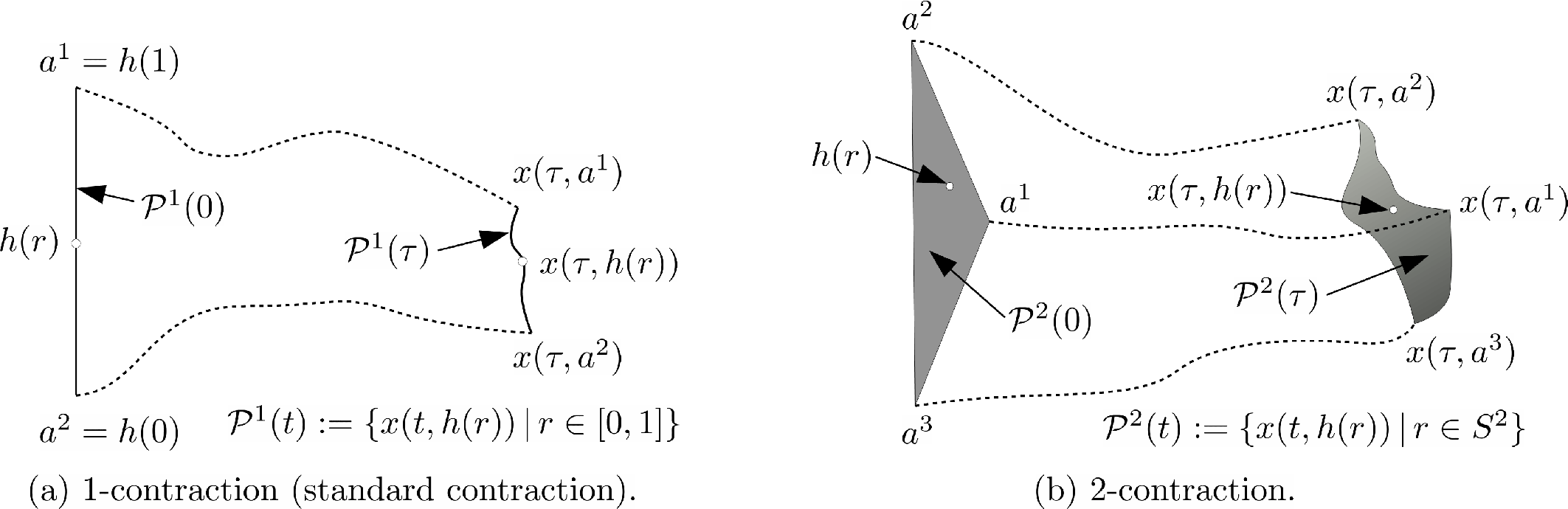}
		\caption{Left: 
	 the length of the curve $\mathcal{P}^1(t)$ decays exponentially (standard contraction).  Right:  in  $2$-contractive 
	 systems 
  the area of the surface~$\mathcal{P}^2(t)$  decays exponentially.}
		\label{illus} 
	\end{center}
\end{figure*}

\subsection{Nonlinear systems}\label{subsec:nonlinear}

Consider the time-varying nonlinear system~\eqref{eq:nonlinsys}.
Pick~$k \in [1,n-1]$. Let
$
S^k:=\{   r\in \R^k \st  r_i\geq 0, r_1+\dots+ r_k\leq 1\}
$
 denote the unit simplex in~$\R^k$. 
Pick~$a^1,\dots,a^{k+1} \in \Omega $. For~$r\in S^k$,   let
$
h(r):=(\sum_{i=1}^k r_i a^i)+(1- \sum_{i=1}^k r_i  )a^{k+1},
$
i.e.  a convex combination of the~$a^i$s, and 
let
\begin{equation}\label{eq:witr}
w^i(t,r):= \frac{\partial   }{\partial r_i} x(t, h(r)) ,\quad i=1,\dots,k.
\end{equation} 
Thus,~$w^i(t,r)$ measures how a change in the initial 
condition~$h(r)$  via a 
change in~$r_i$, affects the solution at time~$t$. 
Note that 
$
w^i(0,r) = \frac{\partial   }{\partial r_i} x(0, h(r)) =a^i-a^{k+1}
$, $i=1,\dots,k$.

\begin{Definition}  \label{def:con_nonli}
The time-varying nonlinear system~\eqref{eq:nonlinsys} is called \updt{\emph{$k$-contractive}}
if there exist~$  \eta>0$ and a vector norm $|\cdot|$  such that 
 for any~$a^1,\dots,a^{k+1} \in \Omega$ and any~$r\in S^k$, \updt{the mapping~$W: \R_+ \times S^k \to \R^{n \times k}$
 defined by
 \[
 W(t,r) := \begin{bmatrix} w^1(t, r) & \cdots &  w^k(t, r) \end{bmatrix}
 \]
 satisfies
	 \be\label{eq:poert}
	 \left| W^{(k)}(t, r)
	\right| \leq   \exp(-\eta t) \left| W^{(k)}(0, r) \right|, \text{ for all } t\geq 0.
	 \ee}
\end{Definition}

To explain  the geometric meaning of this definition, 
pick a domain~$\mathcal{D}\subseteq  S^k$. Then~$k$-contraction implies that
	 \begin{align}\label{eq:kcontnet}
	 &\left| \int_\mathcal{D}   W^{(k)}(t, r)   \diff r  \right|\nonumber \\& \leq  
	 \int_\mathcal{D} \left|  W^{(k)}(t, r) \right| \diff r \nonumber \\&= 
	 \int_\mathcal{D} \left|  \begin{bmatrix}  \frac{\partial }{\partial r_1}x(t,h(r)) & \cdots & \frac{\partial }{\partial r_k}x(t,h(r))   \end{bmatrix}^{(k)}
	 \right| \diff r\\
	 &\leq
	   \exp(-\eta t)
	 \int_\mathcal{D} \left| W^{(k)}(0, r) \right| \diff r \nonumber \\
	 &=  \exp(-\eta t)
	 \left|   \begin{bmatrix} (a^1-a^{k+1})  & \cdots &  (a^k-a^{k+1}) \end{bmatrix}^{(k)}  \right| \int_\mathcal{D} \diff r.\nonumber 
	 \end{align}
	 Note that~\eqref{eq:kcontnet}
	 is the \updt{volume} of the 
	 $k$-surface~$x(t,h(r))$  
	 over the parameter space~$r\in \mathcal{D}$ (see~\eqref{eq:wjskcont}). Thus, \updt{$k$-contraction} implies that this volume  decays to zero at an exponential rate.

\begin{Example}
Suppose that~\eqref{eq:nonlinsys}
is~\updt{$1$-contractive}. 
Pick~$a^1,a^2 \in \Omega$. Then~$
w^1(t,r)=\frac{\partial} {\partial r} x(t, r a^1+(1-r) a^2),
$
and~\eqref{eq:kcontnet} with~$\mathcal{D}=S^1=[0,1]$ becomes
\begin{align*}
 \left| \int_0^1   \frac{\partial} {\partial r} x(t, r a^1+(1-r) a^2)   \diff r  \right| &\leq  
  \exp(-\eta t)
  \left|  a^1-a^{2}  \right |,
\end{align*}
that is,
$
 \left| x(t,a^1)-x(t,a^2)  \right|  \leq  
  \exp(-\eta t)
  \left|  a^1-a^{2} \right |.
$
Thus,~\updt{$1$-contraction} is just contraction. 
\end{Example}

 Fig. \ref{illus}
illustrates   the relation between standard contraction \updt{(i.e. $1$-contraction)} and~\updt{$2$-contraction}.

\begin{Example}
Consider the special case where~$f(t,x)=A(t)x$. Then~\eqref{eq:nonlinsys} is an LTV system. Assume that~$0\in \Omega$. 
Then
\begin{align*}
w^i(t,r)&=\frac{\partial}{\partial r_i} x(t,h(r)) \\
        &=\frac{\partial}{\partial r_i} (\Phi(t)  h(r) )\\
        &=  \Phi(t)  (a^i-a^{k+1}) \\
        &=  x(t,a^i)-x(t,a^{k+1}) ,
\end{align*}
where~$\Phi(\cdot)$ is the transition matrix 
corresponding to the linear dynamics. 
Taking~$a^{k+1}=0$,
Eq.~\eqref{eq:poert} \updt{reduces to   condition~\eqref{eq:contdef} in Definition~\ref{def:conlin}}.
\end{Example}

\subsection{Sufficient conditions for \updt{$k$-contraction} in nonlinear systems}

The next result provides an easy to check sufficient condition for~\updt{$k$-contraction} in terms of the~$k$th additive compound of the Jacobian
of the vector field.
\begin{Theorem} \label{prop:suf_non}
	Suppose that there exist~$\eta>0$ and a vector norm~$|\cdot|$ such that
	\be\label{eq:plost}
	 \mu( J^{[k]}(t,  a ) )\leq -\eta, \text{ for all } a\in\Omega,  t\geq 0.
	\ee
	Then~\eqref{eq:nonlinsys} is \updt{$k$-contractive}.
\end{Theorem}
 
\begin{proof}
	\updt{The definitions of~$W(t,r)$ and~$w^i(t,r)$} give
	 \begin{equation}\label{eq:dwtr}  \begin{aligned}
	 \frac{d}{dt} W(t,r) &= \frac{d}{dt} \frac{\partial x(t,h(r))}{\partial r} \\
	 &= \frac{\partial}{\partial r} f(t,x(t, h(r))) \\
	 &= J(t,x(t,h(r))) \frac{\partial}{\partial r} x(t,h(r))\\
	 &= J(t,x(t, h(r))) W(t,r).
	 \end{aligned}\end{equation}
	 Thus,
	 \[
	 \frac{d}{dt} W^{(k)}(t,r)  = 
	 J^{[k]} (t,x(t, h(r) )) W^{(k)}(t,r),
	 \]
	 and~\eqref{eq:plost}  implies  that \updt{\eqref{eq:poert} holds} for all~$a^1,\dots,a^{k+1} \in \Omega$, $r\in S^k$. 
\end{proof}


If~\eqref{eq:plost} holds for some~$L_p$ norm, with~$p\in\{1,2,\infty\}$, then
arguing as in Corollary~\ref{thm:graded} shows that for any~$\ell \geq k$ we have
$
	 \mu_p( J^{[\ell]}(t,  a ) )\leq -\eta<0$
for all~$t\geq 0$ and all~$a\in \Omega$, so  the nonlinear system is~\updt{$\ell$-contractive}
w.r.t. the~$L_p$ norm.

Recall that~$A\in\R^{n\times n}$ is called Metzler if all its off-diagonal entries are non-negative. 
The nonlinear system~\eqref{eq:nonlinsys}
is called~\emph{$k$-cooperative} if~$J^{[k]}(t,x)$ is Metzler for all~$t \geq 0$ and all~$x \in \Omega $~\citep{Eyal_k_posi}. In other words,~$\dot y=J^{[k]} y$ is a cooperative dynamical system~\citep{hlsmith}. 
Since~$J^{[1]}=J$, this is a generalization of
cooperative  systems (and in fact the 
case~$k=n-1$ corresponds to competitive systems~\citep[Lemma~4]{Eyal_k_posi}).
The next result  provides a sufficient condition for such a system to be~\updt{$k$-contractive} w.r.t. a scaled~$L_1$ norm.   For the special  case~$k=1$, this is closely related to known results on contractive cooperative  systems~\citep{coogan2019contractive}. 
We use~$1_q$ to denote the vector in~$\R^q$ with all entries one.

\begin{Proposition}
Suppose that~\eqref{eq:nonlinsys} is $k$-cooperative. Let~$r := \binom{n}{k}$. 
If there exist~$\eta>0$ and~$v \in \mathbb{R}^r$, with~$v_i>0$ for all~$i$, such that
\be\label{eq:rpoy}
  v^T J^{[k]} (t,x) \leq -\eta 1_r^T \text{ for all } t\geq 0   \text{ and }x\in\Omega
		\ee
	then~\eqref{eq:nonlinsys} 
	is \updt{$k$-contractive} w.r.t. the scaled $L_1$
	norm $\vert x \vert_V: = \vert V x  \vert_1 $, where~$V: = \diag(v_1,   \dots, v_r)$. 
\end{Proposition}
	\begin{proof}
	Let~$ q^T  := 1_r^T V J^{[k]}  V^{-1}$, that is, the entries of the vector~$q$
	are  all the column sums of the 
	matrix~$V J^{[k]}  V^{-1}$.
		Then~$q^T= v^T J^{[k]}  V^{-1}  $, and~\eqref{eq:rpoy}
		implies that~$q^T \leq -1_r^T \eta \min_{i} \{v_i^{-1}\}$. 
	Thus,
	\begin{align*}
	\mu_1(V J^{[k]}  V^{-1}) = \max_{i\in [1,r]} q_i   \leq -\eta \min_{i} \{v_i^{-1}\}
	<0,
	\end{align*}
	where in the first equation we used the fact that since~$J^{[k]} $ is Metzler,
	so is $V J^{[k]}  V^{-1}$. 
\end{proof}

\section{Applications}\label{sec:app}
Li and Muldowney derived several deep 
 results on~$2$-contractive  systems,
although they never used this terminology~\citep{muldo1990,LI2000295,li1993bendixson}.  
These results found many applications in models
 from mathematical epidemiology (see e.g.~\cite{SEIR_LI_MULD1995}). 
These models typically  have at least two equilibrium points,
 corresponding to the disease-free and endemic steady-states. Hence, they cannot 
be~\updt{$1$-contractive} w.r.t. any norm. 
 We begin with  an intuitive presentation of two (somewhat simplified) results that we will use later on, referring to~\citep{muldo1990,li1995} for the full technical details and proofs. 

\subsection{Preliminaries}
The next  result 
is a generalization of  Bendixson's   criterion for the non-existence 
of limit cycles in planar systems.
 
\begin{Theorem} \label{thm:mu2order}\citep{muldo1990}
	Consider the nonlinear time-invariant system:
	\begin{equation} \label{eq:nonli_ti}
	\dot{x} = f(x),
	\end{equation}
	where $f: \mathbb{R}^n \to \mathbb{R}^n$ is $C^1$. Suppose that either
	\begin{equation} \label{eq:2cont}
	\mu\left( J^{[2]}(x) \right)< 0 \text{ for all } x \in \mathbb{R}^n,
		\end{equation}
or
	\begin{equation} \label{eq:2contopp}
\mu\left( -J^{[2]}(x) \right)< 0 \text{ for all } x \in \mathbb{R}^n.
	\end{equation}
 Then \eqref{eq:nonli_ti} has no non-trivial
	periodic solutions.
\end{Theorem}

 Intuitively speaking, the proof is based on the following idea.  Suppose 
 that the system admits a nontrivial periodic solution~$x(t)=x(t+T)$ with minimal period~$T>0$. 
Let~$\gamma$ denote the corresponding invariant curve.
Let~$D$ denote the  trace of a 2-surface whose
boundary is~$\gamma$ and whose surface area is a minimum.
The invariance of~$\gamma$ implies that~$x(t,\gamma)$ 
is also the trace of a 2-surface with boundary~$\gamma$. 
The \updt{$2$-contraction} condition~\eqref{eq:2cont} 
 implies that the area of~$x(t,\gamma)$ is strictly smaller  
 than the area of~$D$ for any~$t>0$.
But this contradicts the definition of~$D$. Condition~\eqref{eq:2contopp}
yields the same contradiction after  replacing~$t$ with~$-t$.

Note that when~$n=2$,
 Condition~\eqref{eq:2cont} [Condition~\eqref{eq:2contopp}]
 becomes~$ \operatorname{div}(f):=\frac{\partial f}{\partial x_1} + \frac{\partial f}{\partial x_2} <0$ [$\operatorname{div}(f)>0$],
so Thm.~\ref{thm:mu2order} is a generalization of Bendixson's theorem for planar systems.

The next result  provides a sufficient condition based on~\updt{$2$-contraction}
guaranteeing 
that an equilibrium   is globally asymptotically  stable.   
\begin{Theorem} \label{thm:wandering} \citep{li1995}
	Consider the nonlinear time-invariant system~\eqref{eq:nonli_ti},
	where $f: \mathbb{R}^n \to \mathbb{R}^n$ is $C^1$.
	Assume that its trajectories evolve on a convex and compact set~$\Omega$, 
	and that
	\[
	\mu \left  (J^{[2]} (x)    \right ) <0 \text{ for all }x \in \Omega. 
	\]
	Then every solution emanating from~$\Omega$ converges to the set of equilibria.
	If in addition there exists a unique equilibrium~$e\in \Omega$ then every solution emanating from~$\Omega$ converges to~$e$.
\end{Theorem}

The proof is based on the following argument. 
Recall that a point~$x_0\in\Omega$ is called \emph{wandering} 
for~\eqref{eq:nonli_ti} if there exists a neighborhood~$U$ of~$x_0$  and  a time~$T>  0 $
such that
\[
					U \cap x(t,U) =\emptyset  \text{ for all } t\geq T. 
\]
In other words, any solution emanating from~$U$     never
returns to~$U$ after time~$T$.
A point~$x_0$ is called \emph{non-wandering} if it is not wandering. Non-wandering points
are important in analyzing the asymptotic behavior of solutions. 
For example, an equilibrium, and  more generally, 
   any point in an omega limit set  is non-wandering. 
Suppose that the conditions in Thm.~\ref{thm:mu2order} hold. 
Assume that~\eqref{eq:nonli_ti} admits a point~$z\in \R^n$ that is non-wandering and is not an equilibrium. 
By the Closing Lemma~\citep{imp_c1}, there exists a~$C^1$ vector field~$\tilde f$, that is arbitrarily close to~$f$ in the~$C^1$ topology,
and~$\dot x=\tilde f(x)$ admits a non-trivial periodic solution.
(Roughly speaking, it is possible to ``close'' the non-wandering trajectory
into a non-trivial periodic trajectory.)
 But, since~$\tilde f$ is arbitrarily close to~$f$   and~$\Omega$ is compact,~$\tilde f$
  also satisfies the \updt{2-contraction} condition in  Thm.~\ref{thm:wandering}
	and thus cannot have a non-trivial periodic solution. 
We conclude that any non-wandering point, and in particular any 
point in an omega limit set,
 must be an equilibrium. 

 The next result   provides a sufficient condition
 for the stability of a non-trivial periodic solution.
\begin{Theorem} \label{thm:suff2}\citep{muldo1990}
	Suppose that  the nonlinear time-invariant 
	system~\eqref{eq:nonli_ti} admits a periodic solution~$\gamma(t)=\gamma(t+T)$
	with minimal period~$T>0$. 
	If the LTV system
 	 	\begin{equation} \label{eq:zfter}
\dot z=    J^{[2]}(\gamma(t))  z 
	\end{equation}
	is asymptotically stable
 then~$\gamma(t)$ is asymptotically orbitally  stable. 

\end{Theorem}
\begin{proof}
By Floquet's theory, the solution of
\be\label{eq:phitr}
\dot \Phi(t)=J(\gamma(t))\Phi(t),\quad \Phi(0)=I,
\ee
can be written as~$\Phi(t)=R(t)\exp(L t)$, where~$R(t)=R(t+T)$ and~$L\in\R^{n\times n}$. 
The eigenvalues~$\lambda_i$, $i=1,\dots,n$,
 of~$L$ are called the characteristic multipliers and one of them,
say~$\lambda_1$, is zero. 
Then
\begin{align}\label{eq:pshytr}
\Phi^{(2)}(t) &= R^{(2)}(t) (\exp(Lt))^{(2)}\nonumber \\
               &=		R^{(2)}(t) 	\exp( L^{[2]}  t)				.	
 \end{align}
where the second equation follows from~\eqref{eq:gkpo}. The eigenvalues of~$L^{[2]}$
are the sum of every pair of eigenvalues of~$L$, and since~$\lambda_1=0$, every~$\lambda_i$, $i=2,\dots,n$,  is an   eigenvalue  of~$L^{[2]}$. 
It follows from~\eqref{eq:phitr} that
\[
\dot \Phi^{(2)}(t)=J^{[2]}(\gamma(t)) \Phi^{(2)}(t), \quad \Phi^{(2)}(0)=I.
\]
The condition in the theorem implies that $\lim_{t\to\infty} \Phi^{(2)}(t) =0$.
Combining this with~\eqref{eq:pshytr} implies that all the eigenvalues of~$L^{[2]}$ have a negative real part,
so in particular, the real part of~$\lambda_i $, $i=2,\dots,n$, is negative.
\end{proof}

Standard contraction can be applied to prove that all trajectories converge to a unique
equilibrium. If a dynamical system admits more than one equilibrium then it is clearly not contractive. 
Yet, it may be~\updt{$k$-contractive}, with~$k>1$, and 
sometimes this can be used to derive a global understanding of the dynamics. 
We demonstrate this using 
  the analysis of a dynamical model
that generalizes the susceptible-exposed-infectious-recovered~(SEIR) 
model studied by~\cite{SEIR_LI_MULD1995}.
Let~$\R^n_+:=\{x\in\R^n \st  x_i\geq 0,\; i=1,\dots,n \}$.

\subsection{Global analysis of a   \updt{3D} system}
Consider the system:
\begin{align}\label{eq:3sus}
					\dot x_1&= -\lambda  f_1(x_1,x_3) +\zeta -\zeta x_1,\nonumber \\
					\dot x_2&=\lambda  f_1(x_1,x_3)  - c x_2 -\zeta x_2 ,\nonumber \\
					\dot x_3&=c x_2 -f_2(x_3) -\zeta x_3,
\end{align}
where  the parameters~$\lambda,\zeta,c$ 
are positive,  and  the   state-space is~$\Omega:=\{x\in\R^3_+ \st  x_1+x_2+x_3\leq 1\}$.
We assume that 
     for any~$x\in \Omega$, we have 
\begin{align}\label{eq:condi_all_x}
f_i(x)   & \geq 0,\; i\in\{1,2\}, \nonumber \\
f_1(x ) &= 0 \text { iff } x_1x_3 = 0 ,  \\
f_2(x_3) &= 0 \text{ iff } x_3 = 0,  \nonumber  
\end{align}
and    for any~$x\in\Int( \Omega)$,  we have
\begin{align}\label{eq:condi_inside_x}
\frac{\partial}{\partial x_j} f_i(x) & >0, \; i\in\{1,2\},\; j\in\{1,2,3\}, \nonumber  \\
\frac{\partial }{\partial x_3}f_1(x_1,x_3) &\leq  \frac{f_1(x_1,x_3) }{ x_3},  \\
 \frac{ f_2(x_3)}{x_3} & \leq \frac{\partial }{\partial x_3} f_2(x_3)  .\nonumber
\end{align}


In the SEIR model,
$f_1(x_1,x_3)=x_1^q x_3^p$, with~$q>0,p\in(0,1]$,  and~$f_2(x_3) = \ell x_3$, with $\ell > 0$, so this indeed  holds. 
\updt{ Note that~$e^1:=\begin{bmatrix}
1&0&0
\end{bmatrix}^T$ is an equilibrium of~\eqref{eq:3sus}. In the SEIR model, this corresponds to the disease-free steady-state. 
The next result analyzes the asymptotic behavior of~\eqref{eq:3sus}.
\begin{Proposition}\label{prop:main}
Suppose that~\eqref{eq:3sus} admits exactly two equilibrium points~$e^1 =\begin{bmatrix}
1&0&0
\end{bmatrix}^T$,
and~$e^2\in \Int(\Omega)$. If~$e^1$ is not 
an omega limit  point for any~$x_0 \in \Int(\Omega)$, and~$e^2$ is   locally asymptotically stable then
\[
\lim_{t\to\infty} x(t,a) =e^2, \text{ for any }  a \in \Int(\Omega).
\]
\end{Proposition}
In other words, the local stability of~$e^2$ implies its global stability in~$\Int(\Omega)$.
In the SEIR model,~$e^2$ corresponds to the endemic steady-state. 
Since~$e^1\in \partial \Omega$, the  property  that~$e^1$ is not 
an omega limit  point for any~$x_0 \in \Int(\Omega)$
can often be verified via conditions guaranteeing that the dynamics is persistent~\citep{SEIR_LI_MULD1995}.
} 

\begin{proof}[Proof of Prop.~\ref{prop:main}]
\updt{The proof consists of several steps: (1)~Using the theory of~$2$-cooperative systems~\citep{Eyal_k_posi}, it is shown   that the system satisfies    the Poincar\'{e}-Bendixson property:
a nonempty compact omega limit set 
 which does not contain any equilibrium points is a closed orbit;
(2)~Using~$2$-contraction   along any non-trivial periodic solution~$\gamma$ it is shown that~$\gamma$ is  asymptotically orbitally stable; (3)~It is shown that~(1) and~(2) imply that the basin of attraction of~$e^2$ includes~$\Int(\Omega)$.
} 

The Jacobian of~\eqref{eq:3sus}
is
\be\label{eq:kjapf}
J=\begin{bmatrix}  
    -\lambda \frac{\partial f_1}{\partial x_1} & 0 &  -\lambda \frac{\partial f_1}{\partial x_3} \\
    \lambda \frac{\partial f_1}{\partial x_1} & -c
		         & \lambda \frac{\partial f_1}{\partial x_3}\\
    0 &   c & -\frac{\partial f_2}{\partial x_3} 
\end{bmatrix} -\zeta I,
\ee
and  Lemma~\ref{lem:poltr} gives
\[
J^{[2]}=\begin{bmatrix}  
    -\lambda \frac{\partial f_1}{\partial x_1} -c  & \lambda \frac{\partial f_1}{\partial x_3}  &  \lambda  \frac{\partial f_1}{\partial x_3} \\
    c &  -\lambda\frac{\partial f_1}{\partial x_1}-\frac{\partial f_2}{\partial x_3} 
		         &   0 \\
    0 & \lambda   \frac{\partial f_1}{\partial x_1} &-c  -\frac{\partial f_2}{\partial x_3}  
\end{bmatrix}-2\zeta I.
\]
Note that~$J^{[2]}(x)$ is Metzler for any~$x \in\Omega$, and irreducible for any~$x \in\Int(\Omega)$. It follows from the results in~\citep{Eyal_k_posi}
 that~\eqref{eq:3sus} is a strongly 2-cooperative system, and thus it 
satisfies 
the Poincar\'{e}-Bendixson property:
a nonempty compact omega limit set 
 which does not contain any equilibrium points is a closed orbit.

Suppose that for some~$x_0\in\Omega$  the omega limit set~$\omega(x_0)$ 
does not contain any equilibrium points. Then~$\omega(x_0)$ 
is a periodic solution~$\gamma (t)$ of~\eqref{eq:3sus}
 with minimal period~$T>0$.  \updt{Our  next goal is to use 2-contraction to  show
 that~$\gamma$ is asymptotically orbitally stable. We require the following result. }
 \begin{Lemma}\label{lem:intgam}
 The periodic solution satisfies~$\gamma(t) \in \Int(\Omega)$ for all~$t\in[0,T)$. 
 \end{Lemma}
 \begin{proof}
We first  show that $\gamma_i(t) \neq 0$ for all~$ i = 1, 2,3$ and~$t \in [0, T)$. 
 Note that for any~$x \in \Omega $ with $x_1 = 0$ we have~$\dot x_1>0$, so~$\gamma_1(t) \neq 0$ for all~$t$. 
 If~$\gamma_2(\tau) = 0$ for some time~$\tau$, then we must have~$\dot{\gamma}_2(\tau) = \lambda f_1(\gamma_1(\tau), \gamma_3(\tau)) \leq 0$. Since  $\gamma_1(\tau) \neq 0$,   $\gamma_3(\tau) = 0$. In this case, the set~$\{x \in \Omega \st x_2 = x_3 = 0  \}$ is forward invariant for~$t \geq \tau$ and~\eqref{eq:3sus} implies that~$\gamma(t)$ converges to the equilibrium point~$e^1$. This contradicts the fact that~$\gamma$ is a non-trivial periodic solution, and thus $\gamma_2(t) \neq 0 $ for all $t$. A similar argument shows that~$\gamma_3(t) \neq 0$ for all $t$. 
Now  suppose that for some~$\tau \in [0,T)$ we have~$\sum_i \gamma_i(\tau)=1$. Then~\eqref{eq:3sus} gives~$\sum_i \dot  \gamma_i(t)=-f_2(\gamma_3(\tau))<0 $. This implies that~$\sum_i \gamma_i(t)<1$ for all~$t\in[0,T)$, and this completes the proof of the lemma.
 \end{proof}

We now show that~$\gamma$ is asymptotically orbitally stable. Consider the system:
\be\label{eq:dynj2}
\dot z(t)=J^{[2]}(\gamma(t))z(t).
\ee
Define~$D(t):=\diag(1,\gamma_2(t)/\gamma_3(t),\gamma_2(t)/\gamma_3(t) )$.
This is well-defined by Lemma~\ref{lem:intgam}.
Let~$p(t):=D(t)z(t)$. 
Then
\begin{align*}
								\dot p&= \dot D z+ D \dot z \\
								&= ( \dot D D^{-1}  + D J^{[2]}(\gamma ) D^{-1} )p  .
\end{align*}
A   calculation gives~$ \dot D D^{-1}=\diag(0,\frac{\dot \gamma_2}{\gamma_2}-\frac{\dot \gamma_3}{\gamma_3},\frac{\dot \gamma_2}{\gamma_2}-\frac{\dot \gamma_3}{\gamma_3} )$.
Note that this implies that~$\int_0^T \dot D (t)D^{-1} (t)\diff t =0$.

Let~$M:=\begin{bmatrix} 1&0&0\\0 &1&1\\ 0&1&-1  \end{bmatrix}$, 
and define   a scaled~$L_\infty$ norm by~$
|y|_{M,\infty}:=|My|_\infty.
$
Then
\begin{align}\label{eq:pt}
								\frac{d^+}{dt^+}  |p|_{M,\infty}& 
								\leq  \mu_\infty (S (\gamma) )|p|_{M,\infty},
\end{align}
with
\begin{align} \label{eq:sgamma}
	S&:= M(  \dot D D^{-1}  + D J^{[2]}  D^{-1} )M^{-1} \nonumber \\
	&=    \dot D D^{-1}  +M  D J^{[2]} D^{-1} M^{-1} \\
	&=					   \dot D D^{-1}  -2\zeta I \nonumber \\
	&+ \begin{bmatrix} 
	-\lambda \frac{\partial f_1}{\partial x_1}-c& 
	\lambda  \frac{\gamma_3}{\gamma_2}\frac{\partial f_1}{\partial x_3} & 0 \\
	\frac{\gamma_2}{\gamma_3}c &-\frac{\partial f_2}{\partial x_3}-\frac{c}{2} &\frac{c}{2}  \\
	\frac{\gamma_2}{\gamma_3} c &-\lambda  \frac{\partial f_1}{\partial x_1}+\frac{c}{2}   &      -\lambda \frac{\partial f_1}{\partial x_1}-\frac{\partial f_2}{\partial x_3} -\frac{c}{2}  \nonumber
		\end{bmatrix}.
\end{align}
\updt{Eq.~\eqref{eq:pt} implies that
\be \label{eq:diniint}
|p(t)|_{M, \infty} \leq \exp \left (\int_0^t \mu_{\infty}(S(\gamma(s))) \diff s\right ) |p(0)|_{M, \infty}.
\ee
}
\updt{By~\eqref{eq:sgamma},} $\mu_\infty(S)=\max\{g_1,g_2\} $, with
\begin{align}\label{eq:g1g2}
g_1&:=-\lambda \frac{\partial f_1}{\partial x_1}-c+\lambda  \frac{\gamma_3}{\gamma_2}\frac{\partial f_1}{\partial x_3}-2\zeta, \nonumber \\
g_2&:=
\frac{\gamma_2}{\gamma_3} c
-\frac{\partial f_2}{\partial x_3}+
\frac{\dot \gamma_2}{\gamma_2}-\frac{\dot \gamma_3}{\gamma_3}-2\zeta.
\end{align}
Using~\eqref{eq:condi_inside_x} gives  
\begin{align*}
g_1& \leq  -c
+\lambda  \frac{f_1}{\gamma_2}-2\zeta.
\end{align*}
The second equation of~\eqref{eq:3sus} gives
$
		\frac{\dot \gamma_2}{\gamma_2 } =\lambda \frac{f_1}{\gamma_2}-c-\zeta,
$
so~$
g_1 
 \leq \frac{\dot \gamma_2}{\gamma_2} 
 - \zeta  .
$
The third equation of~\eqref{eq:3sus} gives
$
		\frac{\dot \gamma_3}{\gamma_3 } =  \frac{c \gamma_2}{\gamma_3}
		-\frac{f_2}{\gamma_3}-\zeta,
$
and combining this with~\updt{\eqref{eq:condi_inside_x} and}~\eqref{eq:g1g2} yields  
$g_2  
 \leq \frac{\dot \gamma_2}{\gamma_2} 
 - \zeta.
$
We conclude that~$\mu_\infty(S) \leq \frac{\dot \gamma_2}{\gamma_2} 
 - \zeta$. \updt{Therefore, 
 \begin{align*}
 \int_0^t \mu_\infty(S(\gamma(s))) \diff s \leq & \log \gamma_2(t) - \log \gamma_2(0) - \zeta t.
 \end{align*}
Since~$\log \gamma_2(t)$ is bounded for all $t \geq 0$, Eq.~\eqref{eq:diniint}} implies that~$\lim_{t\to\infty} p(t)=0$, so~$\lim_{t\to\infty} z(t)=0$. \updt{Since~\eqref{eq:dynj2} is a linear $T$-periodic system, this implies asymptotic stability}, 
and  Thm.~\ref{thm:suff2} implies that~$\gamma$ is asymptotically orbitally stable. 

Summarizing, if for some~$x_0\in\Omega$ we have that~$\omega(x_0)$ 
does not contain any equilibrium points then~$\omega(x_0)$ 
is a non-trivial asymptotically  orbitally  stable periodic solution~$\gamma$ of~\eqref{eq:3sus}.

To complete the proof of Prop.~\ref{prop:main}, let~$B\subset \Omega $ denote the basin of attraction of~$e^2$. Seeking a contradiction, assume that~$\Int(\Omega) \not\subseteq B  $. Then~$
M:=(\partial B ) \cap \Int(\Omega) \not =\emptyset,
$
and~$M$ is an invariant set. Thus, the closure of~$M$ includes a non empty compact omega limit set and the assumptions in the proposition imply that this omega limit set includes no equilibrium points. Thus, it includes a  non-trivial periodic solution~$\gamma$, where~$\gamma$ is in the interior of~$\Omega$ and is asymptotically orbitally stable.
But this is a contradiction, as~$M$ and thus~$\gamma$ is contained in the alpha limit set of~$e^2$. 
This completes the proof of Prop.~\ref{prop:main}.
\end{proof}

\begin{figure}
	\begin{center}
		\centering
	\includegraphics[scale=0.56]{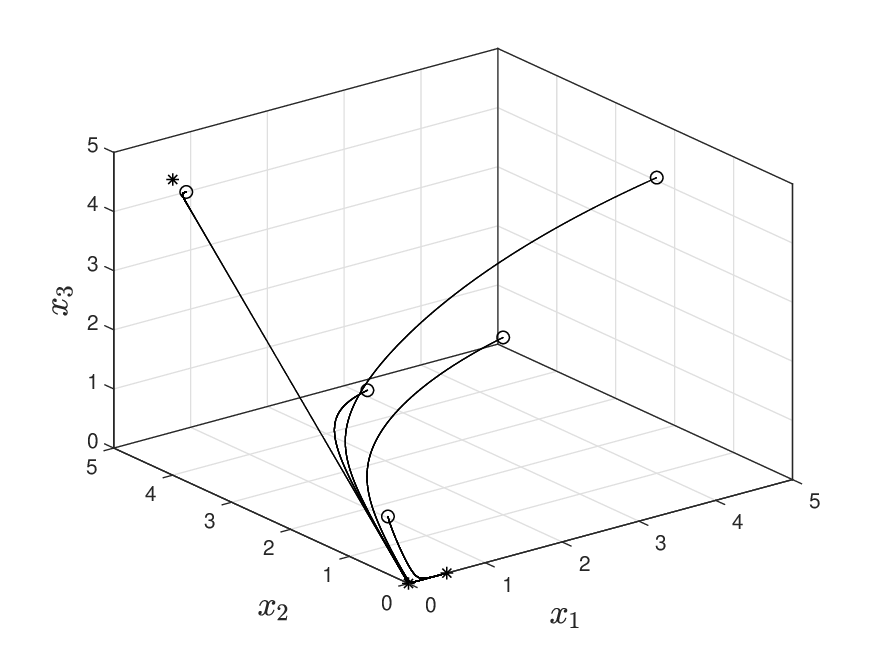}
		\caption{Several trajectories of the Lotka-Volterra system in Example~\ref{exa:lvsys}. 
		The initial conditions [equilibrium points] are marked by $0$ [*]. 
		}
		\label{lvsys} 
	\end{center}
\end{figure}

\updt{
\subsection{$2$-contraction in the 
Lotka-Volterra model}
Consider the Lotka-Volterra model
\be\label{eq:lvol}
\dot x_i=x_i (b_i +\sum_{k=1}^n  a_{ik} x_k),\quad i=1,\dots,n.
\ee
The 
  state-variable~$x_i(t)$    represents
the number of species~$i$ at time~$t$. Note that~$\R^n_+$ is an invariant set of the dynamics. 
This model   has been widely
used in  mathematical 
ecology~\citep{siljak_lrge_scale,sigmund_evolution_book} to study the implications of various interaction patterns
between members of a population sharing a common  habitat.
 
Eq.~\eqref{eq:lvol} can be written   as~$
\dot x=\diag(x_1,\dots,x_n) (b+Ax),
$
where~$b:=\begin{bmatrix}
b_1&\dots & b_n
\end{bmatrix}^T$ and~$A:=\{a_{ij}\}_{i,j=1}^n$.
Thus,~$0$ is an equilibrium,
and if~$A$ is non-singular then~$(-A^{-1}b)$ is an equilibrium.

Let~$\R^n_{++}:=\{x\in\R^n \st x_i>0,\; i\in[1,n]\}$.
There exist known results on the asymptotic behaviour of this model in certain special cases. For example, if~$A$ is diagonally dominant 
(i.e. there exist~$d_i>0$ such that~$d_i a_{ii}+\sum_{j\not =i} d_j  | a_{ji} |<0$ for~$i\in[1,n]$) and there exists an equilibrium~$e\in\R^n_{++}$ then~$\lim_{t\to \infty} x(t,a)=e$ for any~$a\in\Int(\R^n)$ (see e.g.~\cite{Lotka_diag_dominant}).

The quadratic terms in~\eqref{eq:lvol} imply that the model typically admits several equilibrium points and thus cannot be~$1$-contractive.
Our goal is to provide a new sufficient condition 
for $2$-contraction. To do this, 
let~$g_i(x):=b_i+\sum_{k=1}^n a_{ik}x_k$. Then~\eqref{eq:lvol}
 can be written as  the Kolmogorov   system~$\dot x_i=x_i g_i (x)$, $i\in[1,n]$.  
The corresponding Jacobian  
is 
\begin{align*}
J&=\diag(g_1,\dots,g_n) +
\begin{bmatrix}
 x_1\frac{\partial}{\partial x_1} g_1 & \dots&
  x_1\frac{\partial}{\partial x_n} g_1 \\
  &\vdots\\
   x_n\frac{\partial}{\partial x_1} g_n & \dots&
  x_n\frac{\partial}{\partial x_n} g_n
\end{bmatrix} \\
&= 
\diag(g_1,\dots,g_n) +  \diag(x_1,\dots, x_n)A
 , 
\end{align*}
so 
\begin{align}   \label{eq:jac_lv}
J^{[2]}=& \diag(g_1+g_2,g_1+g_3, \dots,g_{n-1}+g_n) 
  \nonumber \\ 
 & +
 (  \diag(x_1,\dots, x_n) A )^{[2]}.
\end{align}
Note that the~$b_i$s   appear only  in  the first matrix on the right-hand side of~\eqref{eq:jac_lv}.
Using this allows to provide easy to verify  sufficient conditions for 2-contraction. Recall that this implies an important asymptotic property, namely, that all bounded solutions converge to an equilibrium. 
The next result demonstrates this for the case~$n=3$.
}
\updt{
\begin{Proposition}\label{prop:lotka}
 Consider~\eqref{eq:lvol} with~$n=3$. If
\begin{align} \label{eq:lvr3con}
& b_i + b_j < 0, \text{ for all } i \neq j, \nonumber \\
& \max\{a_{13} + a_{23},  a_{12} + a_{32}, a_{21}+a_{31}\} \leq 0 , \nonumber \\
& 2 a_{11}+\max \{ a_{21} + |a_{13}|, a_{31} + |a_{12}| \}\leq 0, \\
&2a_{22}+ \max \{ a_{12} + |a_{23}|, a_{32} + |a_{21}| \} \leq 0, \nonumber \\
&  2a_{33}+\max \{ a_{13} + |a_{32}|, a_{23} + |a_{31}| \} \leq 0. \nonumber 
\end{align}
Then   the system is $2$-contractive w.r.t. the~$L_\infty$ norm. 
\end{Proposition}

\begin{proof}
For~$n=3$, Eq.~\eqref{eq:jac_lv} becomes
\begin{align*}
J^{[2]}=\begin{bmatrix}
 m_1&  a_{23}x_2 & -a_{13} x_1\\
 a_{32} x_3 &m_2 & a_{12} x_1\\
 -a_{31} x_3 & a_{21} x_2 & m_3
\end{bmatrix},
\end{align*}
where 
$
      m_1 := b_1 + b_2 +
      (2 a_{11}  + a_{21}) x_1 + (  2 a_{22}+a_{12}  ) x_2 + (a_{13} + a_{23}) x_3$,  
 $
 m_2:=b_1 + b_3 + (2 a_{11} + a_{31}) x_1 + (a_{12} + a_{32}) x_2 + (  2 a_{33}+a_{13} ) x_3$, and 
 $
 m_3:=b_2 + b_3 + (a_{21} + a_{31}) x_1 + (2 a_{22}  + a_{32} )x_2 +(   2 a_{33}+a_{23} ) x_3.
$

Condition~\eqref{eq:lvr3con} ensures that
\[
\mu_\infty  (J^{[2]}(x)) \leq \max\{b_1 + b_2,b_1 + b_3,b_2 + b_3 \} < 0,
\]
for all~$x$,
so the system  is $2$-contractive. 
\end{proof}

\begin{Example} \label{exa:lvsys}
Consider~\eqref{eq:lvol} with~$n=3$, $b_1 = 1$, $b_2 = -2$, $b_3 = -2$, and 
$
A =
\begin{bmatrix}
-2 & -3 & -3 \\
-1 & -2 & 2 \\
1 & 3 & -2
\end{bmatrix}.
$
The dynamics has 
  three equilibrium points in~$\R^3_+$:
  the origin, $\begin{bmatrix}    0&4&5 \end{bmatrix}^T$, and $\begin{bmatrix}    1/2&0&0 \end{bmatrix}^T$.

  Note that~$A$   is not Hurwitz, so it is not diagonally dominant. In fact, the system admits unbounded solutions. 
However, condition~\eqref{eq:lvr3con} holds. Hence, this system is $2$-contractive on~$\R^3_+$, so every bounded trajectory converges to an equilibrium point (see  Fig.~\ref{lvsys}). 
\end{Example}
 }

Our next application of~$k$-contraction 
is to control synthesis.
For a symmetric matrix~$S$, we write~$S\prec  0$ [$S\preceq 0$] if~$S$ is negative   definite [negative semi-definite]. 

\subsection{Control design in a \updt{$2$-contractive} system}
Consider the  affine nonlinear time-invariant     system:
\begin{equation} \label{eq:aff}
\dot{x} = f(x) + G(x)u,
\end{equation}
where $f: \R^n \to \R^n$, $G: \R^n \to \R^{n \times m}$ are $C^1$, and $u \in \R^m$ is the control input.
Let~$J(x) : = \frac{\partial f}{\partial x}(x)$.
We assume that there exists a positive definite matrix $P \in \R^{n \times n}$  such that
\begin{equation} \label{eq:incre}
P J(x) + J^T(x) P \preceq   0 \text{ for all } x \in \R^n.
\end{equation}
 Note that \eqref{eq:incre} is equivalent to $\mu_2(P^{\frac12} J(x) P^{-\frac12}) \leq 0$ for all $x \in \R^n$. 

Define~$V:\R^n\times\R^n \to \R_+$ by
\[
	V(a, b) := \frac12 (a-b)^T P (a-b).
	\]
If~$G(x)\equiv G$, i.e.  the input matrix  is constant,  and~$G$ is full rank 
  then~\eqref{eq:incre} implies that the system \eqref{eq:aff} is \emph{incrementally passive} \citep{van2017L2}
	w.r.t.
	the incremental   storage function 
	$
	V(x(t,a), x(t,b)) $ and the output~$y(x) : = G^T P x$ 	(see also~\cite{pavlov2008incremental, wu2019robust}).
In this case, consider 
the problem of steering the system's output
 to a value~$y(e)$ for some pre-specified~$e\in\R^n$. 
Let $\tilde{x} := x - e$ and $\tilde{y} (x):= y(x) -y(e)$. Then, 
the control design 
\be\label{eq:clu}
u := -kG^T P \tilde{x} + u^*,
\ee
 with $k>0$,
 and $u^*$ satisfying $f(e) + G u^* = 0$ gives
\begin{equation}
\dot{V}(x,e) = \tilde{x}^T P A(x) \tilde{x}^T - k \tilde{y}^T \tilde{y},
\end{equation}
where $A(x):  = \int_0^1 J(e + s \tilde{x} ) \diff s$. 
 Eq. \eqref{eq:incre} yields~$\dot{V}(x,e) \leq - k \tilde{y}^T \tilde{y} \leq 0$. By \emph{LaSalle's invariance principle}, this implies that every solution of the closed-loop system~\eqref{eq:aff} and~\eqref{eq:clu}
converges to the set~$\mathcal{M}$ which  is
 the largest invariant set contained in~$\{x\in \R^n \st \tilde{y}(x) = 0 \}$.
If~$\mathcal{M} =\{e\}$ then~$e$ is GAS.

The next result shows how 
  \updt{$2$-contraction} allows to extend this
	control design method  when the input matrix is allowed to be state-dependent. 

\begin{Proposition}
	Suppose that \eqref{eq:incre}  holds and also  that
	\updt{
	\begin{equation} \label{eq:mu2jx}
	P^{(2)} J^{[2]}(x) + (J^{[2]}(x))^T P^{(2)} \prec  0
	\text{ for all }x\in\R^n,
	\end{equation}}
and that there  exists a $C^1$ mapping $\theta: \R^n \to \R^m$ such that 
	\begin{equation} \label{eq:mugx}
	\mu_2 \left( P^{\frac12} \frac{\partial  }{\partial x}(G(x) \theta(x) )
	P^{-\frac12} \right) \leq 0\text{ for all }x\in\R^n.
	\end{equation} 
 Consider the control $u := \theta(x)$. 
\updt{Then every trajectory of the closed-loop system converges to an equilibrium. If the closed-loop system admits  a unique equilibrium~$e$ then~$e$ is GAS.}
\end{Proposition}

\begin{proof}
	Let $f_c(x): = f(x) + G(x) \theta(x)$, $J_c(x) := \frac{\partial f_c }{\partial x}(x)$, and~$g_c(x):=G(x)\theta(x)$.
Recall that any matrix measure is sub-additive, i.e., $\mu(A+B) \leq \mu(A)+ \mu(B)$ for any~$A,B \in  \R^{n \times n}$ (see e.g.~\cite{feedback2009}), and combining this
	with \eqref{eq:incre} and \eqref{eq:mugx}
		gives
		\begin{align*}
	\mu_2(P^{\frac12} J_c(x) P^{-\frac12}) & \leq \mu_2(P^{\frac12} J(x) P^{-\frac12})\\
	&+ \mu_2 \left(P^{\frac12}  \frac{\partial g_c }{\partial x}    P^{-\frac12} \right)\\& \leq 0.
	\end{align*}
	This implies that the closed-loop system
	is globally uniformly bounded. 
	Hence, for any initial condition~$a \in \R^n$, 
	there exists a compact set~$\mathcal{D}=\mathcal{D}(a)$ such that~$x(t, a) \in \mathcal{D}$ for all $t \geq 0$.
	
	By \eqref{eq:matirxm}, \eqref{eq:matirxm_k}, and \eqref{eq:mugx},
	$
    \mu_2 \left( \left(P^{\frac12}  \frac{\partial g_c } {\partial x}  P^{-\frac12}\right)^{[2]} \right)  
    \leq  \mu_2 \left( P^{\frac12} \frac{\partial g_c } {\partial x}  P^{-\frac12} \right) \leq 0.
	$
	\updt{Since $P$ is positive definite,   so are $P^{(2)}$ and $(P^{(2)})^{\frac12}$. Hence, \eqref{eq:mu2jx} ensures that $\mu_2((P^{(2)})^{\frac12} J^{[2]}  (P^{(2)})^{-\frac12}) < 0$.
	So}
	\begin{align*}
	 \mu_2((P^{\frac12} &J_c  P^{-\frac12} )^{[2]}) \\
	=& \mu_2 \left((P^{\frac12} J  P^{-\frac12} )^{[2]} + \left(P^{\frac12}  \frac{\partial g_c } {\partial x}   P^{-\frac12}\right)^{[2]} \right) \\
	\leq& \mu_2((P^{\frac12} J  P^{-\frac12} )^{[2]}) + \mu_2 \left( \left(P^{\frac12}  \frac{\partial g_c } {\partial x}  P^{-\frac12}\right)^{[2]} \right)\\
	\leq& \mu_2((P^{(2)})^{\frac12} J^{[2]}  (P^{(2)})^{-\frac12}) + \mu_2( P^{\frac12}  \frac{\partial g_c } {\partial x}  P^{-\frac12}) \\
	<& 0.
	\end{align*}
     Thm.~\ref{prop:suf_non} implies that the closed-loop system
	is \updt{$2$-contractive} w.r.t. a scaled $L_2$ norm, and 
 Thm.~\ref{thm:wandering} implies \updt{that every trajectory converges to an equilibrium.}
\end{proof}

The above control design requires solving
the partial differential equation~\eqref{eq:mugx}. 
In certain cases, numerical algorithms can be used to design~$\theta(x)$. For example, if~$G(x)$ is   a polynomial  and we also parameterize~$\theta(x)$ as a polynomial, then sum of squares programming may be  efficient. This approach has been used in the context of~$1$-contraction theory, see e.g.~\citep{aylward2008stability}.
 
\section{Conclusion}

Contraction theory has found numerous applications in systems and control theory. 
However, it is clear that this theory is too restrictive for many systems. For example, if a system admits more than one equilibrium point then it is not contractive  w.r.t. any norm. 

We considered  a geometric  generalization of contraction theory called~\updt{$k$-contraction}.
For the special case~$k=1$ this reduces to standard contraction. 
An easy to check sufficient 
condition  for~\updt{$k$-contraction} is that  
some matrix measure of  
 the~$k$th additive compound of the Jacobian is uniformly negative. 
 \updt{In the case of~$1$-contraction, it is known that under certain regularity conditions the Jacobian  condition is in fact not only sufficient but also   necessary for contraction~\cite[Prop. 3]{sontag_cotraction_tutorial}. An interesting open problem is  whether this condition  is also necessary for $k$-contraction.}
  
We described several implications of \updt{$k$-contraction} to the asymptotic analysis of nonlinear dynamical systems and to control synthesis. 
  To the best of our knowledge, this is the first application of~\updt{$k$-contraction}, with~$k>1$,
in  control theory.  
We believe that \updt{$k$-contraction}, with~$k>1$,
can be used to address  various   system and control problems   
 for dynamical models  where standard contraction  theory  
cannot be applied.

Standard contraction  implies
 \emph{entrainment}
in nonlinear systems with a  time-varying and periodic vector field~\citep{entrain2011,entrain_master}. 
This is important in many applications. For example,
 synchronous generators must entrain to the frequency of the grid.
Biological organisms 
 must develop internal clocks that entrain to the~24h solar day, and so on. 
An important research direction is to  study the implications of~\updt{$k$-contraction}
in dynamical systems with a time-varying and periodic vector field. 

\subsection*{Acknowledgments}
\updt{
We are grateful to the anonymous reviewers and the AE for many helpful comments that helped us to improve the presentation of the results in this paper. We thank J.-J. Slotine for discussions on some of the topics presented here. } 

\section*{Appendix}

\begin{proof}[Proof of Prop.~\ref{prop:sub}]
Let~$\Phi(t)$ be the solution of
$\dot  \Phi(t)=A(t) \Phi(t)$, $\Phi(0)=I_n$. 
Since \eqref{eq:ltv} is uniformly stable, $\Phi(t)$ is uniformly bounded.
Recall that 
$\dot  \Phi^{(k)}(t)=A^{[k]}(t) \Phi^{(k)}(t)$, $\Phi^{(k)}(0)=I_r$, where~$r:=\binom{n}{k}$.

Suppose that Condition~\eqref{cond:nk1} holds. 
Let~$e^i$ be the $i$th canonical vector in~$\R^n$. 
Since~$\dim \mathcal{X} = n-k+1$, there exist~$c_1,\dots,c_k \in \R$, not all zero, 
such that~$ \sum_{i=1}^k c_i e^i\in\mathcal{X}$.
Hence,
\[
0= \lim_{t\to\infty} x(t,\sum_{i=1}^k c_i e^i)=
\lim_{t\to\infty}\sum_{i=1}^k c_i x(t,  e^i).
\]
Combining this  with  the uniform stability assumption  implies that
$
\lim_{t\to \infty } 
\begin{bmatrix}
x(t,e^1)&\dots&x(t,e^k)
\end{bmatrix}^{(k)} =0,
$
that is, 
\[
0=\lim_{t\to \infty }\Phi^{(k)} (t) 
\begin{bmatrix}
e^1&\dots&e^k
\end{bmatrix}^{(k)} .
\]
We conclude that the first column of~$\Phi^{(k)} (t) $ converges to zero. A similar argument shows that this holds for any column of~$\Phi^{(k)} (t) $.    This shows that Condition~\eqref{cond:nk1}  implies Condition~\eqref{cond:sysmkop}.

To prove the converse implication, suppose that 
 Condition~\eqref{cond:sysmkop} holds. 
Pick~$k$   vectors~$a^1,\dots,a^k\in\R^n$.
Define~$X(t) :=\begin{bmatrix} x(t,a^1)& \cdots & x(t,a^k) \end{bmatrix}$. Then  
$
X(t) =\Phi(t) X(0). 
$
By   uniform boundness, there exists an increasing  sequence of times~$t_i$
such that~$\lim_{i\to \infty}t_i=\infty$ and~$P:=\lim_{i\to\infty} X(t_i)$ exists.
Since~$\dot X^{(k)}=A^{[k]}X^{(k)}$,  Condition~\eqref{cond:sysmkop} implies that~$P^{(k)} =0$,
i.e. all minors of order~$k$ of~$P$ are zero. This implies that there exists~$c\in\R^k\setminus\{0\}$ such that
\begin{align*}
0&=Pc \\
&=\lim_{i\to \infty} \sum_{j=1}^k c_j x(t_i,a^j)\\
&= \lim_{i\to \infty}   x(t_i,  \sum_{j=1}^k c_j a^j)\\
&= \lim_{t\to \infty}   x(t,  \sum_{j=1}^k c_j a^j),
\end{align*}
where the last step follows from the uniform stability 
assumption. 
Summarizing, every set of~$k$ linearly independent vectors~$a^1,\dots,a^k \in \R^n$ 
generates a vector~$\sum_{j=1}^k c_j a^j\not =0 $
such that~$\lim_{t\to \infty}   x(t,  \sum_{j=1}^k c_j a^j)=0$.
This  proves that Condition~\eqref{cond:nk1} holds.
\end{proof}


\begin{thebibliography}{48}
\providecommand{\natexlab}[1]{#1}
\providecommand{\url}[1]{\texttt{#1}}
\expandafter\ifx\csname urlstyle\endcsname\relax
  \providecommand{\doi}[1]{doi: #1}\else
  \providecommand{\doi}{doi: \begingroup \urlstyle{rm}\Url}\fi

\bibitem[Achieser(1992)]{achieser1992theory}
N.~I. Achieser.
\newblock \emph{Theory of Approximation (translated by C.J. Hyman)}.
\newblock Dover Publications, Inc., Mineola, New York, 1992.

\bibitem[Aghannan and Rouchon(2003)]{aghannan2003intrinsic}
N.~Aghannan and P.~Rouchon.
\newblock An intrinsic observer for a class of {Lagrangian} systems.
\newblock \emph{IEEE Trans.\ Automat.\ Control}, 48\penalty0 (6):\penalty0
  936--945, 2003.

\bibitem[Aminzare and Sontag(2014)]{sontag_cotraction_tutorial}
Z.~Aminzare and E.~D. Sontag.
\newblock Contraction methods for nonlinear systems: A brief introduction and
  some open problems.
\newblock In \emph{{Proc.\ 53rd IEEE Conf. on Decision and Control}}, pages
  3835--3847, Los Angeles, CA, 2014.

\bibitem[Aylward et~al.(2008)Aylward, Parrilo, and
  Slotine]{aylward2008stability}
E.~M. Aylward, P.~A. Parrilo, and J.-J.~E. Slotine.
\newblock Stability and robustness analysis of nonlinear systems via
  contraction metrics and sos programming.
\newblock \emph{Automatica}, 44\penalty0 (8):\penalty0 2163--2170, 2008.

\bibitem[Coogan(2019)]{coogan2019contractive}
S.~Coogan.
\newblock A contractive approach to separable {Lyapunov} functions for monotone
  systems.
\newblock \emph{Automatica}, 106:\penalty0 349--357, 2019.

\bibitem[Coppel(1965)]{coppel1965stability}
W.~A. Coppel.
\newblock \emph{Stability and {A}symptotic {B}ehavior of {D}ifferential
  {E}quations}.
\newblock Heath, Boston, 1965.

\bibitem[Desoer and Vidyasagar(2009)]{feedback2009}
C.~A. Desoer and M.~Vidyasagar.
\newblock \emph{Feedback Synthesis: Input-Output Properties}.
\newblock SIAM, Philadelphia, 2009.

\bibitem[Do~Carmo(1992)]{carmo1992riemannian}
M.~P. Do~Carmo.
\newblock \emph{Riemannian Geometry}.
\newblock Birkh{\"a}user, 1992.

\bibitem[Fallat and Johnson(2011)]{total_book}
S.~M. Fallat and C.~R. Johnson.
\newblock \emph{Totally Nonnegative Matrices}.
\newblock Princeton University Press, Princeton, NJ, 2011.

\bibitem[Fiedler(2008)]{fiedler_book}
M.~Fiedler.
\newblock \emph{Special Matrices and Their Applications in Numerical
  Mathematics}.
\newblock Dover Publications, Mineola, NY, 2 edition, 2008.

\bibitem[Forni and Sepulchre(2014)]{forni2014}
F.~Forni and R.~Sepulchre.
\newblock A differential {L}yapunov framework for contraction analysis.
\newblock \emph{IEEE Trans.\ Automat.\ Control}, 59\penalty0 (3):\penalty0
  614--628, 2014.

\bibitem[{Forni} and {Sepulchre}(2019)]{forni_diff_diss}
F.~{Forni} and R.~{Sepulchre}.
\newblock Differential dissipativity theory for dominance analysis.
\newblock \emph{IEEE Trans.\ Automat.\ Control}, 64\penalty0 (6):\penalty0
  2340--2351, 2019.

\bibitem[Gantmacher(1960)]{Gantmacher_vol1}
F.~R. Gantmacher.
\newblock \emph{The Theory of Matrices}, volume~I.
\newblock Chelsea Publishing Company, 1960.

\bibitem[Hofbauer and Sigmund(1988)]{sigmund_evolution_book}
J.~Hofbauer and K.~Sigmund.
\newblock \emph{The Theory of Evolution and Dynamical Systems}.
\newblock Cambridge University Press, 1988.

\bibitem[Horn and Johnson(2013)]{matrx_ana}
R.~A. Horn and C.~R. Johnson.
\newblock \emph{Matrix Analysis}.
\newblock Cambridge University Press, 2 edition, 2013.

\bibitem[Jafarpour et~al.(2021)Jafarpour, Cisneros-Velarde, and
  Bullo]{jafarpour2020weak}
S.~Jafarpour, P.~Cisneros-Velarde, and F.~Bullo.
\newblock Weak and semi-contraction for network systems and diffusively-coupled
  oscillators.
\newblock \emph{IEEE Trans.\ Automat.\ Control}, 2021.
\newblock doi: 10.1109/TAC.2021.3073096.

\bibitem[Li and Muldowney(1995{\natexlab{a}})]{SEIR_LI_MULD1995}
M.~Y. Li and J.~S. Muldowney.
\newblock Global stability for the {SEIR} model in epidemiology.
\newblock \emph{Math. Biosciences}, 125\penalty0 (2):\penalty0 155--164,
  1995{\natexlab{a}}.

\bibitem[Li and Muldowney(1995{\natexlab{b}})]{li1995}
M.~Y. Li and J.~S. Muldowney.
\newblock On {R. A. Smith's} autonomous convergence theorem.
\newblock \emph{Rocky Mountain J. Math.}, 25\penalty0 (1):\penalty0 365--378,
  1995{\natexlab{b}}.

\bibitem[Li and Muldowney(1996)]{li1996geometric}
M.~Y. Li and J.~S. Muldowney.
\newblock A geometric approach to global-stability problems.
\newblock \emph{SIAM J. Math. Anal.}, 27\penalty0 (4):\penalty0 1070--1083,
  1996.

\bibitem[Li and Muldowney(2000)]{LI2000295}
M.~Y. Li and J.~S. Muldowney.
\newblock Dynamics of differential equations on invariant manifolds.
\newblock \emph{J. Diff. Eqns.}, 168\penalty0 (2):\penalty0 295--320, 2000.

\bibitem[Li and Muldowney(1993)]{li1993bendixson}
Y.~Li and J.~S. Muldowney.
\newblock On {Bendixson's} criterion.
\newblock \emph{J. Diff. Eqns.}, 106:\penalty0 27--39, 1993.

\bibitem[Lohmiller and Slotine(1998)]{LOHMILLER1998683}
W.~Lohmiller and J.-J.~E. Slotine.
\newblock On contraction analysis for non-linear systems.
\newblock \emph{Automatica}, 34:\penalty0 683--696, 1998.

\bibitem[Lohmiller and Slotine(2000)]{doi:10.1002/aic.690460317}
W.~Lohmiller and J.-J.~E. Slotine.
\newblock Nonlinear process control using contraction theory.
\newblock \emph{AIChE Journal}, 46\penalty0 (3):\penalty0 588--596, 2000.

\bibitem[Lu(1998)]{Lotka_diag_dominant}
Z.~Lu.
\newblock Global stability for a {Lotka-Volterra} system with a weakly
  diagonally dominant matrix.
\newblock \emph{Appl. Math. Lett.}, 11\penalty0 (2):\penalty0 81--84, 1998.

\bibitem[Manchester et~al.(2018)Manchester, Tang, and
  Slotine]{manchester2018unifying}
I.~R. Manchester, J.~Z. Tang, and J.-J.~E. Slotine.
\newblock Unifying robot trajectory tracking with control contraction metrics.
\newblock In A.~Bicchi and W.~Burgard, editors, \emph{Robotics Research: Volume
  2}, pages 403--418. Springer International Publishing, 2018.

\bibitem[Margaliot and Sontag(2019)]{margaliot2019revisiting}
M.~Margaliot and E.~D. Sontag.
\newblock Revisiting totally positive differential systems: A tutorial and new
  results.
\newblock \emph{Automatica}, 101:\penalty0 1--14, 2019.

\bibitem[Margaliot et~al.(2014)Margaliot, Sontag, and Tuller]{RFM_entrain}
M.~Margaliot, E.~D. Sontag, and T.~Tuller.
\newblock Entrainment to periodic initiation and transition rates in a
  computational model for gene translation.
\newblock \emph{PLoS ONE}, 9\penalty0 (5):\penalty0 e96039, 2014.

\bibitem[Margaliot et~al.(2016)Margaliot, Sontag, and
  Tuller]{3gen_cont_automatica}
M.~Margaliot, E.~D. Sontag, and T.~Tuller.
\newblock Contraction after small transients.
\newblock \emph{Automatica}, 67:\penalty0 178--184, 2016.

\bibitem[Margaliot et~al.(2017)Margaliot, Tuller, and Sontag]{cast_book}
M.~Margaliot, T.~Tuller, and E.~D. Sontag.
\newblock Checkable conditions for contraction after small transients in time
  and amplitude.
\newblock In N.~Petit, editor, \emph{Feedback Stabilization of Controlled
  Dynamical Systems: In Honor of {Laurent} {Praly}}, pages 279--305. Springer
  International Publishing, Cham, Switzerland, 2017.

\bibitem[Margaliot et~al.(2018)Margaliot, Gr{\"u}ne, and
  Kriecherbauer]{entrain_master}
M.~Margaliot, L.~Gr{\"u}ne, and T.~Kriecherbauer.
\newblock Entrainment in the master equation.
\newblock \emph{Royal Society Open Science}, 5\penalty0 (4):\penalty0 172157,
  2018.

\bibitem[Muldowney(1990)]{muldo1990}
J.~S. Muldowney.
\newblock Compound matrices and ordinary differential equations.
\newblock \emph{The Rocky Mountain J. Math.}, 20\penalty0 (4):\penalty0
  857--872, 1990.

\bibitem[Muldowney(1998)]{mol_appl_comp}
J.~S. Muldowney.
\newblock Compound matrices and applications.
\newblock 1998.
\newblock URL
  \url{https://www.researchgate.net/publication/326273499_Compound_Matrices_and_Applications}.
\newblock Lecture notes for Universidad de Los Andes, Merida, Venezuela.

\bibitem[Pavlov and Marconi(2008)]{pavlov2008incremental}
A.~Pavlov and L.~Marconi.
\newblock Incremental passivity and output regulation.
\newblock \emph{Systems \& Control Letters}, 57\penalty0 (5):\penalty0
  400--409, 2008.

\bibitem[Pugh(1967)]{imp_c1}
C.~C. Pugh.
\newblock An improved closing lemma and a general density theorem.
\newblock \emph{American J. Math.}, 89\penalty0 (4):\penalty0 1010--1021, 1967.

\bibitem[Russo and di~Bernardo(2009)]{russo2009solving}
G.~Russo and M.~di~Bernardo.
\newblock Solving the rendezvous problem for multi-agent systems using
  contraction theory.
\newblock In \emph{Proc. of the $48$h IEEE Conference on Decision and Control
  held jointly with $28$th Chinese Control Conference}, pages 5821--5826. IEEE,
  2009.

\bibitem[Russo et~al.(2010)Russo, di~Bernardo, and Sontag]{entrain2011}
G.~Russo, M.~di~Bernardo, and E.~D. Sontag.
\newblock Global entrainment of transcriptional systems to periodic inputs.
\newblock \emph{PLOS Computational Biology}, 6:\penalty0 e1000739, 2010.

\bibitem[Sanfelice and Praly(2011)]{sanfelice2011convergence}
R.~G. Sanfelice and L.~Praly.
\newblock Convergence of nonlinear observers on {$\mathbb{R}^n$} with a
  {R}iemannian metric (part {I}).
\newblock \emph{IEEE Trans.\ Automat.\ Control}, 57\penalty0 (7):\penalty0
  1709--1722, 2011.

\bibitem[Schwarz(1970)]{schwarz1970}
B.~Schwarz.
\newblock Totally positive differential systems.
\newblock \emph{Pacific J. Math.}, 32\penalty0 (1):\penalty0 203--229, 1970.

\bibitem[Siljak(2007)]{siljak_lrge_scale}
D.~D. Siljak.
\newblock \emph{Large-Scale Dynamic Systems: Stability and Structure}.
\newblock Dover Publications, 2007.

\bibitem[Slotine and Wang(2005)]{slotine2005study}
J.-J.~E. Slotine and W.~Wang.
\newblock A study of synchronization and group cooperation using partial
  contraction theory.
\newblock In V.~Kumar, N.~Leonard, and A.~S. Morse, editors, \emph{Cooperative
  Control}, volume 309 of \emph{Lecture Notes in Control and Information
  Science}, pages 207--228. Springer, Berlin, Heidelberg, 2005.

\bibitem[Smith(1995)]{hlsmith}
H.~L. Smith.
\newblock \emph{Monotone Dynamical Systems: An Introduction to the Theory of
  Competitive and Cooperative Systems}, volume~41 of \emph{Mathematical Surveys
  and Monographs}.
\newblock Amer. Math. Soc., Providence, RI, 1995.

\bibitem[Smith(1986)]{smith_Hausdorff_dim}
R.~A. Smith.
\newblock Some applications of {Hausdorff} dimension inequalities for ordinary
  differential equations.
\newblock \emph{Proc. Royal Society of Edinburgh: Section A Mathematics},
  104:\penalty0 235--259, 1986.

\bibitem[Strom(1975)]{storm1975}
T.~Strom.
\newblock On logarithmic norms.
\newblock \emph{SIAM J. Numerical Analysis}, 12\penalty0 (5):\penalty0
  741--753, 1975.

\bibitem[Teschl(2012)]{teschl2012ordinary}
G.~Teschl.
\newblock \emph{Ordinary Differential Equations and Dynamical Systems}.
\newblock American Mathematical Soc., 2012.

\bibitem[van~der Schaft(2017)]{van2017L2}
A.~van~der Schaft.
\newblock \emph{$L_2$-gain and passivity techniques in nonlinear control}.
\newblock London, Springer, third edition, 2017.

\bibitem[Vidyasagar(1978)]{vid}
M.~Vidyasagar.
\newblock \emph{Nonlinear Systems Analysis}.
\newblock Prentice Hall, Englewood Cliffs, NJ, 1978.

\bibitem[Weiss and Margaliot(2021)]{Eyal_k_posi}
E.~Weiss and M.~Margaliot.
\newblock A generalization of linear positive systems with applications to
  nonlinear systems: Invariant sets and the {Poincar\'{e}-Bendixson} property.
\newblock \emph{Automatica}, 123:\penalty0 109358, 2021.

\bibitem[Wu et~al.(2019)Wu, van~der Schaft, and Chen]{wu2019robust}
C.~Wu, A.~van~der Schaft, and J.~Chen.
\newblock Robust trajectory tracking for incrementally passive nonlinear
  systems.
\newblock \emph{Automatica}, 107:\penalty0 595--599, 2019.

\end{thebibliography}
 \end{document}